\documentclass[11pt]{article}
\usepackage[dvips]{color}
\usepackage{hyperref}
\usepackage{mathrsfs}
\usepackage{graphicx}
\usepackage{subfigure}
\usepackage{amsfonts}
\usepackage{amsmath}
\usepackage{amssymb}
\usepackage{booktabs}
 \usepackage{epstopdf}
\definecolor{dgreen}{rgb}{0,.8,.3}
\definecolor{blue}{rgb}{.2,.3,.7}
\definecolor{red}{rgb}{1.0,0.2,0.2}
\usepackage{color}
\input amssym.def
\usepackage{multirow}
\usepackage[pagewise]{lineno}
\usepackage{tikz}
\usetikzlibrary{arrows,shapes,chains} 

\input epsf.tex

\begin{document}

\renewcommand{\Box}{\rule{2.2mm}{2.2mm}}
\newcommand{\BOX}{\hfill \Box}

\newtheorem{eg}{Example}[section]
\newtheorem{thm}{Theorem}[section]
\newtheorem{lemma}{Lemma}[section]
\newtheorem{example}{Example}[section]
\newtheorem{remark}{Remark}[section]
\newtheorem{proposition}{Proposition}[section]
\newtheorem{corollary}{Corollary}[section]
\newtheorem{defn}{Definition}[section]
\newtheorem{alg}{Algorithm}[section]
\newtheorem{ass}{Assumption}[section]
\newenvironment{case}
    {\left\{\def\arraystretch{1.2}\hskip-\arraycolsep \array{l@{\quad}l}}
    {\endarray\hskip-\arraycolsep\right.}

\def\argmin{\mathop{\rm argmin}}

\makeatletter
\renewcommand{\theequation}{\thesection.\arabic{equation}}
\@addtoreset{equation}{section} \makeatother

\title{A Surrogate Value Function Formulation for Bilevel Optimization\thanks{\baselineskip 9pt This paper is supported by National Key  R\&D Program of China (Nos. 2021YFA1000300, 2021YFA1000301), the National Natural
Science Foundation of China (Nos. 12071342, 12021001, 92473208, 12071108) and the Natural Science Foundation of Hebei Province (No. A2025202027).}}
\author{Mengwei Xu\thanks{\baselineskip 9pt Institute of Mathematics, Hebei University of Technology, Tianjin 300401, China. 
E-mail: xumengw@hotmail.com},\ 
Yu-Hong Dai\thanks{\baselineskip 9pt Academy of Mathematics and Systems Science, Chinese Academy of Sciences, Beijing 100190, China. E-mail: dyh@lsec.cc.ac.cn},\   
 Xin-Wei Liu\thanks{\baselineskip 9pt Institute of Mathematics, Hebei University of Technology, Tianjin 300401, China. E-mail: mathlxw@hebut.edu.cn} and\  
 Meiqi Ma\thanks{\baselineskip 9pt Institute of Mathematics, Hebei University of Technology, Tianjin 300401, China. E-mail:  meiqima@126.com}
}  

\date{}
\maketitle
{\bf Abstract.} 
The value function formulation captures the hierarchical nature of bilevel optimization through the optimal value function of the lower level problem, yet its implicit and nonsmooth characteristics pose significant analytical and computational difficulties. We introduce a surrogate value function formulation that replaces the intractable value function with an explicit surrogate derived from lower level stationarity conditions. This surrogate formulation preserves the essential idea of the classical value function model but fundamentally departs from Karush–Kuhn–Tucker (KKT) formulations, which embed lower level stationary points into the upper level feasible region and obscure the hierarchical dependence.
Instead, it enforces the hierarchy through a dominance constraint that remains valid even when lower level constraint qualifications fail at the solution. We establish equivalence with the original bilevel problem, reveal the failure of standard constraint qualifications, and show that its strong stationarity implies that of KKT models. To handle the complementarity constraints in the surrogate formulation, we apply a smoothing barrier augmented Lagrangian method and prove its convergence to solutions and Clarke stationary points. Extensive experiments demonstrate the robustness and high numerical precision of this formulation, especially in nonconvex settings, including the classical Mirrlees’ problem where KKT models fail.



\newpage

\baselineskip 18pt
\parskip 2pt
\section{Introduction} 
Bilevel optimization generalizes classical optimization through a hierarchical structure, where the feasible region of the upper level problem is determined by the optimal solutions of a parameterized lower level problem. 
This hierarchical dependency brings substantial computational challenges, as the interaction between the two levels makes the feasible region highly nonconvex.
Bilevel programs naturally arise in diverse areas such as economics, engineering design,  and machine learning \cite{gttzz,KK,ygy,YeYZ}.
Mathematically, a bilevel program (BP) takes the form
\begin{eqnarray*}
({\rm BP})~~~~~~\min_{x,y} && F(x,y)\nonumber\\
{\rm s.t.} 
&& y\in S(x), 
\end{eqnarray*}
where  $S(x)$ denotes the solution set of the lower level program
\begin{eqnarray*}
({\rm P}_x)~~~~~~~~
\min_y && f(x,y)\\
{\rm s.t.} && g(x,y)\leq 0.
\end{eqnarray*}
Here we assume that 
$F:\mathbb{R}^d\times \mathbb{R}^l \rightarrow \mathbb{R}$ 
is continuously differentiable,  
$f:\mathbb{R}^d\times \mathbb{R}^l \rightarrow \mathbb{R}$, $g:\mathbb{R}^d\times \mathbb{R}^l \rightarrow \mathbb{R}^m$
are continuously differentiable and twice continuously differentiable with respect to (w.r.t.) the variable $y$.
To focus on the main ideas, we omit possible constraints on the upper level variable and equality constraints in the lower level, as the analysis can be extended to these cases without much difficulty.

To describe this hierarchical feasible region, various single-level formulations have been developed. Most notably, the value function formulation (VP) provides a clear and rigorous way to express the bilevel structure \cite{fzz,jmz,Outrata, YeYZ}.
It characterizes upper level feasibility through the optimal value of the lower level problem:
\begin{eqnarray}
~~~~~~~~~~~~({\rm VP})~~~~~~\min_{x,y} && F(x,y)\nonumber\\
{\rm s.t.} && f(x,y)-V(x) \leq 0,\label{VP}\\
&& g (x,y) \leq 0,\nonumber
\end{eqnarray}
where $V(x) :=\inf_{y\in Y(x)} \ f(x,y)$ denotes the lower level value function, and $Y(x):=\{y\in\mathbb{R}^l: g(x,y)\le 0\}$ denotes the feasible set of the lower level problem.
This formulation is equivalent to the original bilevel problem, without requiring convexity or constraint qualifications for  the lower level problem \cite{dz}.
Although theoretically rigorous, VP remains difficult to handle in practice.
The value function $V(x)$ is typically implicit, nonsmooth, and nonconvex, and hence difficult to evaluate or compute \cite{yz2}.
The conditions guaranteeing Lipschitz continuity or the existence of subdifferentials are inherently restrictive [31, 41] (see also Proposition \ref{partialV}).
Moreover,  the nonsmooth optimality conditions of VP may fail at locally optimal solutions, since generalized Mangasarian–Fromovitz constraint qualification (MFCQ) is violated at any feasible point \cite[Proposition~3.2]{yz}.
In summary, VP provides an equivalent model of bilevel optimization, but the implicit and nonsmooth value function makes it difficult to analyze, and the optimality conditions may fail in practice.

In addition to the value function approach, another widely studied formulation replaces the lower level problem with its KKT conditions, transforming (BP) into a mathematical program with equilibrium constraints (MPEC) \cite{b98,dz}:
\begin{eqnarray}
({\rm KP})~~~~~~\min_{x,y,s} && F(x,y)\nonumber\\
{\rm s.t.} 
&&   \nabla_{2} f(x,y)+\sum_{i=1}^m s_i \nabla_{2} g_i(x,y)=0, \label{KP}\\
&&s_i\geq 0,\ g_i(x,y)\leq 0,\ s_i g_i(x,y)= 0,\ i=1,\cdots,m, \nonumber
\end{eqnarray}
where 
$\nabla_{2} f$ denotes the gradient of $f$ w.r.t. $y$. 
The KKT reformulation provides an explicit formulation grounded in well-established stationarity theory and a variety of algorithmic frameworks, including regularization, relaxation, and decomposition methods \cite{ass, flp,hks,ips,kdb}. 
However, its validity relies on the convexity and the satisfaction of certain constraint qualifications of the lower level problem \cite{as}.
When these conditions fail at the optimum,  the KKT system may not hold, causing the model to exclude true solutions.
Even in the convex case, local equivalence can be lost when the lower level multipliers are not unique \cite{DempeDu}.

Beyond these two classical formulations, hybrid or duality-based extensions were proposed.
The combined program (CP) \cite{yz2} integrates  KKT and value function constraints, while duality-based methods \cite{llzz,llz,oa} are constructed by Lagrange, Wolfe or Mond-Weir dual representations.
These formulations improve computational tractability and offer new theoretical viewpoints for specific classes of bilevel problems.
Nevertheless, they still face intrinsic limitations and can yield spurious stationary points \cite{dm24}. 
Further developments for structured bilevel programs can be found in \cite{dkpk,dz,gttzz,hss22,hxlt,lm,nwy,ss17} and references therein.

To overcome the analytical difficulties of the value function formulation while preserving its hierarchical dependence on the lower level optimal response, we propose a new model that makes this dependence explicit and computationally tractable.
Our key insight is that for any  \(x\in\mathbb{R}^d\) and any  \(u\in S(x)\), the equality \(V(x)=f(x,u)\) holds.
This observation inspires the Surrogate Value Function (SVF) formulation, which replaces the implicit value function $V(x)$ with the explicit surrogate $f(x,u)$, where  $u$ is a reference stationary point:
\begin{eqnarray}
({\rm SVF})~~~~~~~\min_{x,y,u,s}  && F(x,y)\nonumber\\
{\rm s.t.} && f(x,y)-f(x,u)\leq 0,\label{fl}\nonumber\\
&&\nabla_{2} f(x,u)+  \sum_{i=1}^m s_i \nabla_{2} g_i(x,u) =0,\label{fl1}\\
&& s_i\geq 0,\ g_i(x,u)\leq 0,\ s_i g_i(x,u)=0,\ i=1,\cdots,m,\label{fl2}\\
&& g(x,y)\leq 0.\nonumber
\end{eqnarray}
Here \(s\) denotes the associated multiplier. 
The constraint that $f(x,y)\le f(x,u)$ explicitly encodes the optimality condition inherent in the value function formulation, thereby preserving its structure within a tractable framework.
Although SVF also involves complementarity conditions, its foundation differs fundamentally from KP: rather than enforcing stationarity at the solution $y$, it verifies it at a reference point $u$, thus avoiding the dependence on constraint qualifications at the solution.

As illustrated in Example~\ref{Counterexample}, when the lower level problem violates the KKT conditions, the true global solution lies outside the feasible region of KP, and VP fails to satisfy its optimality conditions, leading to the breakdown of nonsmooth algorithms.
Furthermore, since the regularity conditions in Proposition~\ref{partialV} are violated, the subdifferential of the value function becomes difficult to compute.
In contrast, SVF recovers the correct solution by introducing a representative stationary point $u$ ensuring both feasibility and stationarity.
This example demonstrates that SVF preserves the spirit of VP while remaining valid in degenerate or nonconvex settings where both KP and VP fail.

We establish theoretical properties of SVF that clarify its equivalence and stationarity relations with existing formulations.
First, when the lower level problem is pseudoconvex, SVF is  equivalent to the original bilevel problem.
Next, under the joint convexity of the lower level problem and suitable constraint qualifications, we prove that any Mordukhovich stationary point of SVF corresponds to a stationary point of VP.
Moreover, under these assumptions, a stationary point of VP also induces a strong stationary point of SVF.
If in addition, the lower level KKT conditions are satisfied at the solution, the strong stationarity of SVF implies that of  KP.
Finally, we show that MPEC-type constraint qualifications, such as the MPEC no nonzero abnormal multiplier constraint qualification (MPEC-NNAMCQ), fail to hold at all feasible points of SVF.

To efficiently handle the complementarity structure of SVF, we employ a smoothing barrier augmented Lagrangian (SBAL) method \cite{ld}, which unifies interior-point \cite{bhn,cgr,wb} and augmented Lagrangian \cite{bm,cgt,gpsy} techniques.
This approach avoids maintaining strictly interior iterates and demonstrated strong performance in convex quadratic and semidefinite programs \cite{zld23,zld24}.
The SBAL framework is particularly well suited to SVF, as its complementarity conditions can be smoothly approximated.
We show that solutions and stationary points of the smoothed problems  converge to solutions and Clarke stationary points of SVF as the smoothing parameter tends to zero.
Extensive numerical experiments on 109 benchmark problems confirm the robustness and global convergence of SVF, especially in  nonconvex cases.
Notably, in the classical Mirrlees’ example \cite{yz2}, a well-known counterexample demonstrating the failure of KKT-based models, SVF produces a near-optimal solution with high numerical precision.

The main contributions of this paper are summarized as follows:
\begin{itemize}
 \item[{\rm (i)}]  We propose the SVF formulation, which preserves the philosophy of the value function model while avoiding the intractability of the implicit function and the limitations of KKT formulations.
 
  \item[{\rm (ii)}]  We establish equivalence between SVF and the original bilevel program, prove the failure of standard constraint qualifications, and derive strong stationarity conditions that are strictly stronger than those of KP.
  
\item[{\rm (iii)}]  We apply a SBAL scheme and prove its convergence to solutions and Clarke stationary points of SVF. Extensive experiments demonstrate its robustness and strong performance in nonconvex and degenerate cases where traditional models fail.

\end{itemize}

The remainder of the paper is organized as follows. Section 2 presents foundational concepts and a motivating example. 
Section 3 establishes equivalence and stationarity properties for SVF and compares its stationary points with those of existing reformulations. Section 4 develops the SBAL approximation and analyzes its convergence.
Section 5 reports numerical experiments and comparative results.
Section 6 concludes the paper.

We adopt the following standard notation in this paper. For any two vectors $a$ and $b $ in $\mathbb{R}^n$, we denote by $a^T b$  their inner product. Given a function $G: \mathbb{R}^n\rightarrow \mathbb{R}^m$, we denote its Jacobian by $\nabla G(z)\in \mathbb{R}^{m\times n}$ and, if $m=1$, the gradient $\nabla G(z)\in \mathbb{R}^n$ is considered as a column vector.   
For a set $\Omega\subseteq \mathbb{R}^n$, we denote by  
co $\Omega$
the convex hull of $\Omega$.
For a matrix $A\in \mathbb{R}^{n\times m}$, $A^T$ denotes its transpose. 
Let $e_i\in \mathbb{R}^m$ be the vector such that the $i$th component is 1 and others are 0.

\section{Preliminaries and An Illustrative Example}
This section presents essential preliminaries. We first review the differentiability of
 the value function  and then define stationarity notions for SVF. 
 Finally, we introduce an illustrative example that reveals  limitations of existing formulations and motivates the proposed SVF formulation.
 
\subsection{Stationarity Concepts}
For  a Lipschitz continuous function $f: \mathbb{R}^d \rightarrow  \mathbb{R}$,
 the Bouligand subdifferential (B-subdifferential) of $f$ at $x$ is denoted by
 $\partial f(x) := \{ v \in \mathbb{R}^d \mid \exists x_k\in D_f, x_k \to x, \nabla f(x_k) \to v \text{ as } k \to \infty \}$, where $D_f$ denotes the set of points where $f$ is differentiable.
The  Clarke subdifferential (C-subdifferential) of $f$ at ${x}$ is
$ \partial^c f({x}):={co} \partial f({x}).$

The next result summarizes the Lipschitz continuity and subdifferential characterization of the value function $V(x)$ of the bilevel optimization problem (BP).


 \begin{proposition}\label{partialV} \cite{m1,yw}
(i) If $S$ is inner semicontinuous at $x^*$ for some $y^*\in S(x^*)$, i.e., for every $x_k\to x^*$, there is a sequence $y_k\in S(x_k)$  converging to $y^*$ and the MFCQ holds at $y^*$,
then the value function $V(x)$ is Lipschitz continuous near $x^*$ and 
\begin{eqnarray*}
\partial^c V(x^*) &\subseteq & co W_{y^*}(x^*),
\end{eqnarray*}
where $W_{y^*}(x^*):=  \left \{\nabla_{1} f(x^*, y^*)+s^T \nabla_{1} g(x^*,y^*): s\in M(x^*, y^*) \right \}$ and
\begin{eqnarray*}
M(x^*,y)&:=& \left\{ s: 
\begin{array}{c}
0=\nabla_{2} f(x^*,y)+s\nabla_{2} g(x^*,y)\\
s\geq 0,\quad  \langle g(x^*,y),s\rangle=0
\end{array}
\right\}.
\end{eqnarray*}
(ii) Assume that $Y(x)$ is uniformly bounded around $x^*$, i.e., there exists a neighborhood $\mathbb{U}$ of $x^*$ such that the set $\bigcup_{x\in \mathbb{U}} Y(x)$ is bounded. Suppose that MFCQ holds at $y$ for all $y\in S(x^*)$. Then the value function $V(x)$ is Lipschitz continuous near $x^*$ and 
$$\partial^c V(x^*) \subseteq co \bigcup_{y\in S(x^*)} W_{y^*}(x^*).$$
(iii) If $f(x,y)$ is jointly convex, $g(x,y)$ is jointly quasiconvex and the Abadie constraint qualification holds for $Y(x)$ at $(x^*,y^*)$, where $y^*\in S(x^*)$, then $V(x)$ is a convex function and
\begin{eqnarray*}
\partial^c V(x^*) = co W_{y^*}(x^*).
\end{eqnarray*}
\end{proposition}

In the rest of this subsection, we introduce the stationary points for the  problem SVF, which follows from the definitions of MPEC stationary points \cite{y05}. Denote by $\cal{F}_S$ the feasible region of SVF.
  {Given}  $(\bar x,\bar y,\bar u,\bar s)\in \cal{F}_S$, define the index sets:
 \begin{eqnarray*}
I(\bar x,\bar y):=\{i: g_i(\bar x,\bar y)=0\},&&
 I_g(\bar x,\bar u):=\{i: g_i(\bar x,\bar u)=0, \bar{s}_i >0\},\\
 I_0(\bar x,\bar u):=\{i: g_i(\bar x,\bar u)=0, \bar{s}_i =0\},&&
 I_s(\bar x,\bar u):=\{i: g_i(\bar x,\bar u)<0, \bar{s}_i=0\}.
\end{eqnarray*}

\begin{defn}\label{stationarypoint}
We say $(\bar x,\bar y,\bar u,\bar s)\in\cal{F}_{S}$ is a weak stationary point (W-stationary point) of the problem SVF if there exist  $\lambda\in\mathbb{R}_+^{m+1}$, $\mu^{\phi}\in \mathbb{R}^l, \mu^g\in \mathbb{R}^m$  such that $\lambda_0(f(\bar x, \bar y)-f(\bar x, \bar u))=0$,
\begin{eqnarray}
&&0= \nabla_{1} F(\bar x, \bar y)  
+\lambda_0(\nabla_{1} f(\bar x, \bar y)-\nabla_{1} f(\bar x, \bar u) ) +\sum_{i\in I(\bar x,\bar y)} \lambda_i \nabla_{1} g_i(\bar x, \bar y) \nonumber\\
&&~~~+ [\nabla_{21}^2 f + \sum_{i\in I(\bar x,\bar u)} \bar{s}_i \nabla_{21}^2 g_i
](\bar x, \bar u)^T\mu_{\phi}
+ \sum_{i\in I(\bar x,\bar u)} \mu_i^g \nabla_{1} g_i(\bar x, \bar u), \label{lkt-xb} \\
&&0=
\nabla_{2} F(\bar x, \bar y) +\lambda_0 \nabla_{2} f(\bar x, \bar y)+\sum_{i\in I(\bar x,\bar y)} \lambda_i \nabla_{2} g_i(\bar x, \bar y), \label{lkt-yb} \\
&&0= 
-\lambda_0 \nabla_{2} f(\bar x, \bar u)
+\sum_{i\in I(\bar x,\bar u)} \mu_i^g \nabla_{2} g_i(\bar x, \bar u)+  [\nabla_{22}^2 f + \sum_{i\in I(\bar x,\bar u)} \bar{s}_i \nabla_{22}^2 g_i](\bar x, \bar u)^T\mu_{\phi},\label{lkt-ub}\\
&&\mu_i^g=0, i\in I_s(\bar x,\bar u),\quad
\nabla_{2} g_i(\bar x, \bar u)^T\mu_{\phi}=0,\ i\in I_g(\bar x,\bar u). 
\label{lkt-sb1}
\end{eqnarray}
We say $(\bar x,\bar y,\bar u,\bar s)\in\cal{F}_{S}$ is a Clarke stationary point (C-stationary point) of the problem SVF if there exist  $\lambda\in\mathbb{R}_+^{m+1}$, $\mu^{\phi}\in \mathbb{R}^l, \mu^{\psi}\in \mathbb{R}^m$  such that (\ref{lkt-xb})-(\ref{lkt-sb1}) and the following condition hold:
\begin{eqnarray}
&&\mu_i^g\nabla_{2} g_i(\bar x, \bar u)^T\mu_{\phi}\geq 0,\ i\in I_0(\bar x,\bar u).\label{lkt-sbc}
\end{eqnarray}
We say $(\bar x,\bar y,\bar u,\bar s)\in\cal{F}_{S}$ is a Mordukhovich stationary point (M-stationary point) of  SVF if there exist  $\lambda\in\mathbb{R}_+^{m+1}$, $\mu^{\phi}\in \mathbb{R}^l, \mu^{\psi}\in \mathbb{R}^m$  such that (\ref{lkt-xb})-(\ref{lkt-sb1}) and the following condition hold:
\begin{eqnarray}
&&\mu_i^g>0,\ \nabla_{2} g_i(\bar x, \bar u)^T\mu_{\phi}>0,\ {\rm or}\ \mu_i^g\nabla_{2} g_i(\bar x, \bar u)^T\mu_{\phi}= 0,\ i\in I_0(\bar x,\bar u).\label{lkt-sbm}
\end{eqnarray}
We say $(\bar x,\bar y,\bar u,\bar s)\in\cal{F}_{S}$ is a strong stationary point (S-stationary point)  if there exist  $\lambda\in\mathbb{R}_+^{m+1}$, $\mu^{\phi}\in \mathbb{R}^l, \mu^{\psi}\in \mathbb{R}^m$  such that (\ref{lkt-xb})-(\ref{lkt-sb1}) and the following condition hold:
\begin{eqnarray}
&&\mu_i^g\geq 0,\ \nabla_{2} g_i(\bar x, \bar u)^T\mu_{\phi}\geq 0,\  i\in I_0(\bar x,\bar u).\label{lkt-sbs}
\end{eqnarray}
\end{defn}

\subsection{An illustrative Counterexample} 
To demonstrate the failure of classical bilevel reformulations in certain settings, we present a simple example where both value function and KKT approaches fail to recover the global solution.

\begin{eg}\label{Counterexample} 
	\begin{figure}
		\centering
		\includegraphics[width=0.6\textwidth]{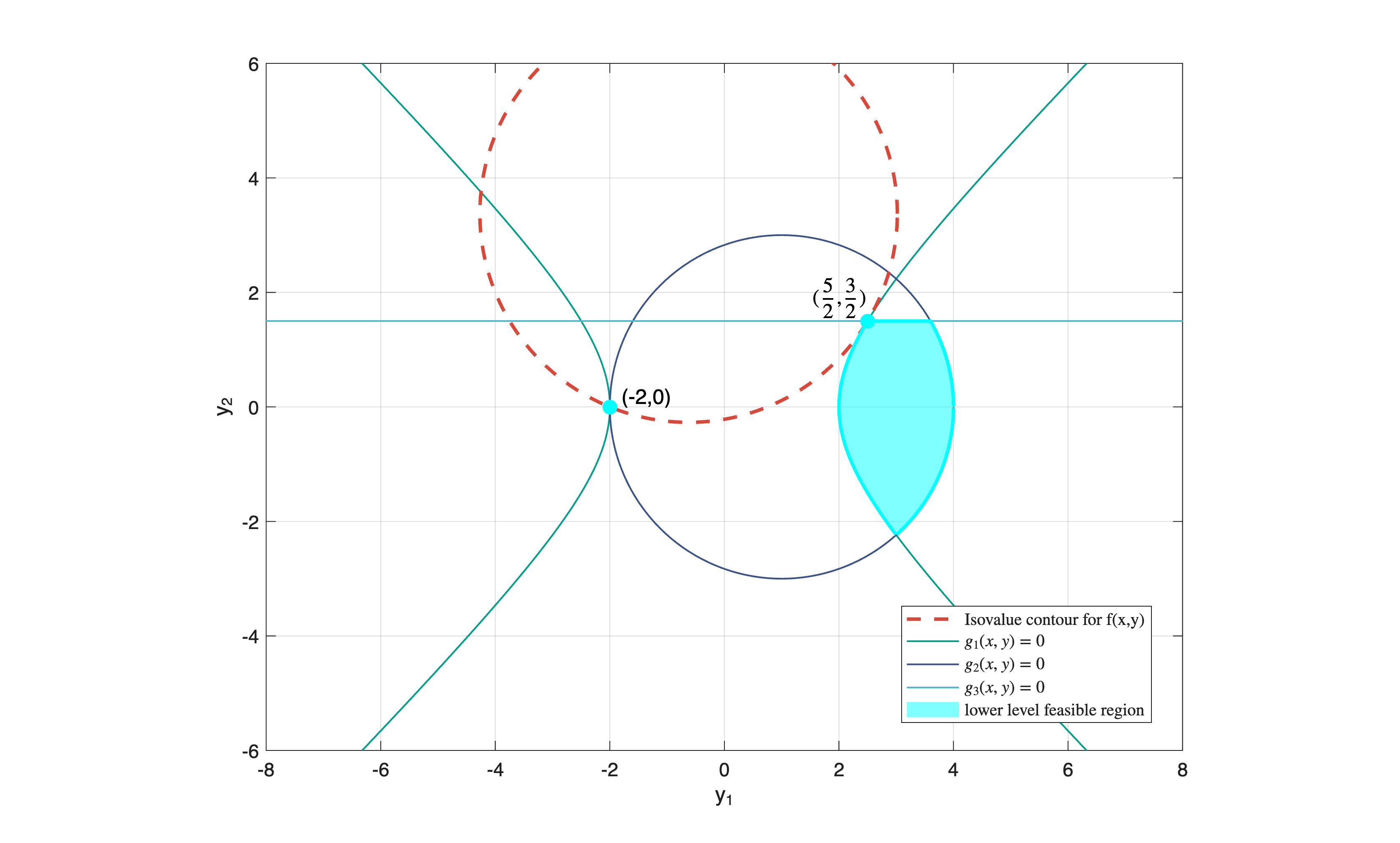}
		\caption{ Feasible region of the lower level problem in  Example \ref{Counterexample} at $x=0$}
		\label{fig:e1}
	\end{figure}
	Consider the bilevel problem:
	\begin{align*}
		\min_{x,y}~~&F(x,y):=x^2 -2x+ 3y_1+y_2 \\
		\rm{s.t.}~~&G_1(x):=x-1\leq 0,\ G_2(x):=-x-1\leq 0,\\
		&y\in \arg\min_y \{f(x,y):= \frac{1}{2}[(y_1-x+\frac{5}{8})^2+(y_2-\frac{27}{8})^2]: \\
		&~~~~~~~~~~~~~~~~
		g_1(x,y):=-{y_1^2}+{y_2^2}+4\leq 0,\\
		&~~~~~~~~~~~~~~~~g_2(x,y):=(y_1-1)^2+y_2^2-9\le0,\\
		&~~~~~~~~~~~~~~~~g_3(x,y):=x+y_2-\frac{3}{2}\le0
		\}.
	\end{align*}
	As shown in Figure \ref{Counterexample}, when $x = 0$, the feasible region of the lower level problem consists of the shaded area together with an isolated feasible point $\bar y:=(-2,0)$. 
	Therefore,  the solution set of the lower level problem at $x=0$ is $S(0)=\{\bar y, \bar u\}$, where $\bar u: =(\frac{5}{2},\frac{3}{2})$. For $-1\le x<0$, the solution set of the lower level problem is the singleton $S(x) = \{\bar y\}$. For $0<x\le 1$,  the unique solution $(y_1(x), y_2(x))$  of the lower level problem is given by the intersection of the right branch of the hyperbola $g_1(x,y)=0$ with the line $g_3(x,y)=0$, thus $y_1(x)=\sqrt{y_2(x)^2+4},\ y_2(x)=\frac{3}{2}-x$. 	By comparing the corresponding upper level objective values across all cases, it follows that the global optimal solution of Example \ref{Counterexample} is $(\bar{x},\bar{y})$ with $\bar x:=0$ and $\bar y=(-2,0)$.

We first consider the problem VP. From the above discussion, it is easy to see that
\begin{eqnarray*}
	V(x)=\begin{cases}
		\frac{1}{2}[(\bar y_1-x+\frac{5}{8})^2+(\bar y_2-\frac{27}{8})^2]  &\text{ if } x<0, \\
		\frac{1}{2}[(y_1(x)-x+\frac{5}{8})^2+(y_2(x)-\frac{27}{8})^2]&\text{ if } x\ge0.
	\end{cases}
\end{eqnarray*}
Since the regularity conditions in Proposition~\ref{partialV} are violated, we compute the C-subdifferential of the value function via the convex hull of its B-subdifferential \cite[Theorem 2.5.1]{c} and derive that $\partial^c V(\bar x)\in 
co \{-\frac{25}{8}, \frac{11}{8}\}$. It is straightforward to compute that   $I(\bar x,\bar y)=\{1,2\}$, $\nabla F(\bar{x},\bar{y})=(-2,3,1)^T$ and 
\begin{eqnarray*}
 \nabla_2 f(\bar x,\bar y)=\binom{-\frac{11}{8}}{-\frac{27}{8}},\ \nabla_2 g_1(\bar x,\bar y)=\binom{4}{0},\  \nabla_2 g_2(\bar x,\bar y)=\binom{-6}{0}.
\end{eqnarray*}		
It is concluded that there exists no  $\lambda_i\geq 0$, $i=0,1,2$ such that the KKT conditions for VP hold:
\begin{eqnarray*}
	&&0\in \nabla_{1} F(\bar x, \bar y)  +\lambda_0(\nabla_{1} f(\bar x, \bar y)-\partial^c V(\bar x)), \label{vp1}\\
	&&0=\nabla_{2} F(\bar x, \bar y) +\lambda_0 \nabla_{2} f(\bar x, \bar y)+\lambda_1 \nabla_{2} g_1(\bar x, \bar y)+\lambda_2 \nabla_{2} g_2(\bar x, \bar y)\label{vp2}.
\end{eqnarray*}
As a result, the global solution $(\bar x,\bar y)$ is not a stationary point of $VP$ and thus the value function formulation cannot solve this problem, as we highlighted in Section 1.	

	We can verify that KKT conditions fail at $\bar{y}$ for the problem $({\rm P}_{\bar x})$. 
	Obviously, there exists no  $\lambda_i\geq 0$, $i=1,2$ such that
	\begin{eqnarray*}
		\nabla_2f(\bar{x},\bar{y})+\lambda_1\nabla_2g_1(\bar{x},\bar{y})+\lambda_2\nabla_2g_2(\bar{x},\bar{y})=0.
	\end{eqnarray*}
  Thus  in this example, $(\bar{x},\bar{y})$ is not  a feasible point of KP and  this global solution of (BP) will be missed by KP.
	
	We now reformulate the original problem into its SVF model:
	\begin{eqnarray*}
		\min_{x,y,u,s}  && F(x,y)\\
		{\rm s.t.} &&g_0(x,y,u):=\frac{1}{2}[(y_1-x+\frac{5}{8})^2+(y_2-\frac{27}{8})^2]-\frac{1}{2}[(u_1-x+\frac{5}{8})^2+(u_2-\frac{27}{8})^2]\leq 0,\\
		&&h_1(x,u,s)=u_1-x+\frac{5}{8}-2 s_1 u_1 +2s_2(u_1-1)=0,\\ 
		&&h_2(x,u,s)=u_2-\frac{27}{8} +2 s_1 u_2 + 2s_2 u_2+s_3=0,\\
		&& s_i\geq 0,\ g_i(x,u)\leq 0,\ s_i g_i(x,u)=0,\ i=1,2,3,\\
		&& g(x,y)\leq 0,\ G(x)\leq0.\nonumber
	\end{eqnarray*}
It is straightforward to verify that the lower level KKT condition holds at $\bar{u}$ with multiplier  $(\bar{s}_1,\bar{s}_2,\bar{s}_3)$, where $\bar{s}_1=\frac{5}{8}, \bar{s}_2=\bar{s}_3=0$ followed from easily calculation:
\begin{eqnarray*}
	\nabla_2 f(\bar x,\bar u)=\binom{\frac{25}{8}}{-\frac{15}{8}},\ \nabla_2 g_1(\bar x,\bar u)=\binom{-5}{3},\ \nabla_2 g_3(\bar x,\bar u)=\binom{0}{1}.
\end{eqnarray*}
Therefore, $(\bar x,\bar y, \bar u,\bar s)$ is a feasible point  of the SVF model.
	By selecting $\lambda :=(\frac{8}{27},0,\frac{35}{81})^T$, $\mu_g:=(-\frac{221}{1242},0,\frac{12}{23})^T$, $\mu_{\phi}:=(-\frac{10}{69},-\frac{50}{207})^T$ and $\lambda_G=(0,0)^T$, 
	we can verify that $(\bar x,\bar y, \bar u, \bar s)$ satisfies the optimality conditions (\ref{lkt-xb})–(\ref{lkt-sb1}) with $I(\bar x,\bar y) =\{1,2\}$, $I_g(\bar x,\bar u)=\{ 1\}$, $I_0(\bar{x},\bar{u})=\{3\}$, and $I_s(\bar x,\bar u)=\{2\}$. Hence, $(\bar{x}, \bar{y}, \bar{u}, \bar{s})$ is a W-stationary point of the SVF model.

						\end{eg}  	
						
This example underscores the fundamental merit of the SVF formulation.
By introducing a reference stationary point and separating the roles of $y$ and $u$, SVF preserves the hierarchical structure of the bilevel problem while overcoming the difficulties caused by VP and KP.
Hence, SVF provides a computationally tractable framework that remains applicable to degenerate  bilevel problems.

\section{Analysis of the SVF Formulation}
This section investigates the theoretical properties of the SVF formulation. 
In Section 3.1, we first prove the  equivalence between SVF and the original bilevel problem under pseudoconvexity.
We then analyze  the failure of MPEC-NNAMCQ, and employ the partial calmness condition to derive M‑stationarity conditions for SVF.
The relationships between the stationarity conditions of SVF and those of existing formulations will be discussed in Section 3.2.

\subsection{Optimality and Stationarity Analysis }
Denote by $L(x,u,s):=f(x,u)+  \sum_{i=1}^m s_i  g_i(x,u)$ the Lagrangian function of the lower level problem. 
We give the following assumption, which holds automatically if the lower level problem is convex.

\begin{ass}\label{pseudoconvex}
Assume that the Lagrangian function $L(x,u,s)$ is pseudoconvex for each $x$. In other words, for any $u_1, u_2\in Y(x)$  and $s\in M(x,u_2)$, $(u_1-u_2)^T \nabla_{2} L(x,u_2,s)\geq 0$ implies $L(x,u_1,s)\geq L(x,u_2,s)$.
\end{ass}



\begin{thm}\label{glosol}
Assume Assumption \ref{pseudoconvex} holds.
Let $(\bar x,\bar y)$  be a local (global) solution of  (BP). 
Suppose $(\bar u,\bar s)$ is a KKT pair of $(P_{\bar x})$. 
Then  $(\bar x,\bar y,\bar u,\bar s)$ is a local (global) solution of $SVF$.
Conversely, suppose $(\bar x,\bar y,\bar u,\bar s)$ is a solution of $SVF$ on $U(\bar x,\bar y)\times  \mathbb{R}^{l+m}$ and for any $(x,y)\in U(\bar x,\bar y)\cap \cal{F}$, there exists a KKT pair of $(P_x)$, where $U(\bar x,\bar y)$ is the neighborhood of $(\bar x,\bar y)$.
Then $(\bar x,\bar y)$  is a local optimal solution of  (BP).
\end{thm}
 {\bf Proof.}  For any $(x,y,u,s)\in \cal{F}_S$, we have that $\nabla_{2} L(x,u,s)=0$.  It follows from the pseudoconvexity of $L(x,\cdot,s)$ that for $y_x\in S(x)$,
  \begin{eqnarray}\label{vl}
V(x)=f(x,y_x)\geq L(x,y_x,s)\geq L(x,u,s)=f(x,u)\geq f(x,y).
\end{eqnarray}
Then $u, y\in S(x)$ and thus $(x,u), (x,y)\in \mathcal{F}$, where $\mathcal{F}$ denotes the feasible region of (BP).

Let $(\bar x,\bar y)$  be a local optimal solution of  (BP), i.e.,
  \begin{eqnarray}\label{solcon1}
F(x,y)\geq F(\bar x,\bar y),\ \forall (x,y)\in U(\bar x,\bar y)\cap \cal{F}.
\end{eqnarray}
Since $(\bar u,\bar s)$ is a KKT pair of $(P_{\bar x})$, we have that  
 \begin{eqnarray*}
&&\nabla_{2} f(\bar x,\bar u)+  \sum_{i=1}^m \bar{s}_i \nabla_{2} g_i(\bar x,\bar u) =0,\\
&& \bar{s}_i\geq 0,\ g_i(\bar x,\bar u)\leq 0,\ \bar{s}_i g_i(\bar x,\bar u)=0,\ i=1,\cdots,m.
\end{eqnarray*}
Then $(\bar x,\bar y,\bar u,\bar s)$ is a feasible point of SVF. 

For any $(x,y,u,s)\in \cal{F}_S$ such that $(x,y)\in U(\bar x,\bar y)$, it follows from (\ref{vl}) that
$y\in S(x)$ and thus $(x,y)$ is a feasible point of $(BP)$.  From (\ref{solcon1}),  $(\bar x,\bar y,\bar u,\bar s)$ is a local solution of $SVF$.

Conversely, suppose that $(\bar x,\bar y,\bar u,\bar s)$ is a solution of $SVF$ on $U(\bar x,\bar y)\times \mathbb{R}^{l+m}$, i.e.,
  \begin{eqnarray}\label{solcon2}
F(x,y)\geq F(\bar x,\bar y),\ \forall (x,y,u,s)\in U(\bar x,\bar y)\times \mathbb{R}^{l+m}\cap \cal{F}_S.
\end{eqnarray}
From (\ref{vl}), $(\bar x,\bar y)\in\cal{F}$.
For any $(x,y)\in U(\bar x,\bar y)\cap \cal{F}$, from statement, there exists a KKT pair $(u,s)$ of $(P_x)$ and thus $(x,y,u,s)\in \cal{F}_S$. From (\ref{solcon2}), $(\bar x,\bar y)$  is a local optimal solution of  (BP).
 \BOX

A major challenge in analyzing the formulation VP is that standard constraint qualifications, such as nonsmooth MFCQ  fail at every feasible point \cite{yz,yz2}. The next theorem shows that our SVF formulation shares the same structural limitation.
Denote by $\mathcal{F}'_{EC}:=\{z: g(z) \leq 0,\  h(z) =0,\ 0\leq G(z)\perp H(z) \geq 0\}$, where
$g:\mathbb{R}^n\rightarrow \mathbb{R}^p, h: \mathbb{R}^n \rightarrow \mathbb{R}^q$, $G, H:\mathbb{R}^n\rightarrow \mathbb{R}^m$ are continuously differentiable.
We now recall definition of MPEC-NNAMCQ.

\begin{defn}[MPEC-NNAMCQ]\cite{y05}
We say that the 
{\rm MPEC-NNAMCQ}  holds at $z^* \in\mathcal{F}'_{EC}$,
if there is no nonzero vector $\lambda=(\lambda^g,\lambda^h, \lambda^G,\lambda^H) \in \mathbb{R}^{p+q+2m}$ such that 
 \begin{eqnarray*}
&&0\in \sum_{i\in I_g} \lambda_i^g \nabla  g_i(z^*) +\sum_{i=1}^q \lambda_i^h \nabla h_i(z^*) -\sum_{i=1}^m[ \lambda_i^G \nabla G_i(z^*) +\lambda_i^H \nabla H_i(z^*)],\\
&& \lambda_i \geq 0 \ \  i \in I_g, \quad \lambda_i^G=0 \ \  i\in \gamma,\quad  \lambda_i^H=0 \ \  i \in \alpha, \\
&&  \mbox{either } \lambda_i^G>0, \lambda_i^H>0 \mbox{ or } \lambda_i^G\lambda^H_i=0  \ \ \forall i \in \beta, \nonumber
\end{eqnarray*}
where $I_g:=I_g(z^*)=\{i: g_i(z^*)=0\},\ \alpha:=\alpha(z^*)=\{i: G_i(z^*)=0, H_i(z^*) >0\},\
  \beta:=\beta(z^*)=\{i: G_i(z^*)=0, H_i(z^*) =0\},\ \gamma:=\gamma(z^*)=\{i: G_i(z^*)>0, H_i(z^*) =0\}$.
\end{defn}

\begin{thm}\label{nna-fail}
The MPEC-NNAMCQ for the problem SVF is violated at any feasible point if the Assumption \ref{pseudoconvex} holds.
\end{thm}
{\bf Proof.}  Denote by $$\widetilde{\cal{F}}:=\{(x,y,u,s):\nabla_{2} f(x,u)+s^T \nabla_{2} g(x,u)=0,\ 0\geq g(x,u)\perp s\geq 0,\ g(x,y)\leq 0\}.$$ 
Suppose $(\bar x,\bar y,\bar u,\bar s)\in \cal{F}_{S}$. Assume that MPEC-NNAMCQ holds for the system $\widetilde{\cal{F}}$ at $(\bar x,\bar y,\bar u,\bar s)$, 
otherwise MPEC-NNAMCQ for the problem SVF is violate automatically.

Since the Assumption \ref{pseudoconvex} holds, similarly as (\ref{vl}), for any $(x,y,u,s)\in \widetilde{\cal{F}}$, we have that $u\in S(x)$ and thus $f(x,y)\geq f(x,u)$.

Since $(\bar x,\bar y,\bar u,\bar s)\in \cal{F}_S$, $f(\bar x,\bar y)-f(\bar x,\bar u)=0$. Then $(\bar x,\bar y,\bar u,\bar s)$ is a minimizer of the problem
\begin{eqnarray*}
\min_{(x,y,u,s)\in \widetilde{\cal{F}}}  & f(x,y)-f(x,u )
\end{eqnarray*}
From the MPEC-NNAMCQ for the system $\widetilde{\cal{F}}$, the M-stationary condition holds at $(\bar x,\bar y,\bar u,\bar s)$ \cite[Corollary 2.1]{y05},
 which means that 
 MPEC-NNAMCQ fails for SVF.
 \BOX

 The above theorems show that when the lower level problem is pseudoconvex, SVF is equivalent to the original bilevel program, yet any feasible point fails to satisfy MPEC tailored constraint qualifications such as MPEC-NNAMCQ.
The following result shows that if the partial calmness condition holds for VP, then it also holds for SVF.
%

\begin{thm}\label{calmvp}
Assume the Assumption \ref{pseudoconvex} holds. Suppose that $(\bar x,\bar y)$ is a local  solution of  VP
and $(\bar x,\bar y,\bar u,\bar s)$ is a local solution of $SVF$.  
If  VP is partially calm at $(\bar x,\bar y)$, then SVF is partially calm at $(\bar x,\bar y,\bar u,\bar s)$, i.e., there exists $\lambda>0$ such that $(\bar x,\bar y,\bar u,\bar s)$ is a local optimal solution of
 \begin{eqnarray*}
({\rm SVF_{\lambda}})~~~~~~~\min_{(x,y,u,s)\in \widetilde{\cal{F}}}  && F(x,y)+\lambda (f(x,y)-f(x,u)),\nonumber
\end{eqnarray*}
where $\widetilde{\cal{F}}$ is denoted in Theorem \ref{nna-fail}.
\end{thm}
{\bf Proof.}  Assume the partially calmness of  VP holds at $(\bar x,\bar y)$, then there exists $\bar{\lambda}>0$ such that for any $\lambda>\bar{\lambda}$, $(\bar x,\bar y)$  is a local optimal solution of 
$$({\rm VP_{\lambda}})~~~~~\min_{x,y}~F(x,y)+\lambda (f(x,y)-V(x))\quad {\rm s.t.}~g(x,y)\leq 0.$$
For any $(x,y,u,s)\in U(\bar x,\bar y,\bar u,\bar s)\cap \widetilde{\cal{F}}$, it follows from Assumption \ref{pseudoconvex} that $u\in S(x)$ and
 \begin{eqnarray*}
 F(\bar x,\bar y)+\lambda (f(\bar x,\bar y)-f(\bar x,\bar u))
&=&F(\bar x,\bar y)+\lambda (f(\bar x,\bar y)-V(\bar x))\\
&\leq& F(x,y)+\lambda (f(x,y)-V(x))\\
&=& F(x,y)+\lambda (f(x,y)-f(x,u)).
\end{eqnarray*}
Then $(\bar x,\bar y,\bar u,\bar s)$ is a local optimal solution of $({\rm SVF_{\lambda}})$ and 
we complete the proof.
\BOX

We now derive stationarity conditions for SVF.  To focus on the main contribution, we state only a sufficient condition for M‑stationarity.
\begin{thm}\label{opt}
Suppose that $(\bar x,\bar y,\bar u,\bar s)$ is a local solution of SVF satisfied the partially calmness condition. If MPEC-NNAMCQ holds for the system $\widetilde{\cal{F}}$ at $(\bar x,\bar y,\bar u,\bar s)$, then it is an M-stationary point of $SVF$.
\end{thm}
{\bf Proof.} From the partially calmness of SVF, we assume there exists $\widetilde\lambda>0$ such that $(\bar x,\bar y,\bar u,\bar s)$ is a local solution of $({\rm SVF_{\widetilde\lambda}})$.
From \cite[Corollary 2.1]{y05}, 
 the Fritz-John type M-stationary condition holds at $(\bar x,\bar y,\bar u,\bar s)$. Thus there exist $r\geq 0$ $\lambda\in\mathbb{R}_+^{m}$, $\mu^{\phi}\in \mathbb{R}^l, \mu^g, \mu^s\in \mathbb{R}^m$ not all equal to zero such that
\begin{eqnarray}
&&0=  r\left[\nabla_{1} F(\bar x, \bar y)+
\widetilde\lambda(\nabla_{1} f(\bar x, \bar y)-\nabla_{1} f(\bar x, \bar u) ) \right]
 +\sum_{i\in I(\bar x,\bar y)} \lambda_i \nabla_{1} g_i(\bar x, \bar y) \label{kt-xb}\nonumber\\
&&~~~+[\nabla_{21}^2 f + \sum_{i\in I(\bar x,\bar u)} \bar{s}_i \nabla_{21}^2 g_i
](\bar x, \bar u)^T\mu_{\phi}
+ \sum_{i\in I(\bar x,\bar u)} \mu_i^g \nabla_{1} g_i(\bar x, \bar u), \nonumber \\
&&0=r\left[\nabla_{2} F(\bar x, \bar y)+
\widetilde\lambda \nabla_{2} f(\bar x, \bar y)\right]+\sum_{i\in I(\bar x,\bar y)} \lambda_i \nabla_{2} g_i(\bar x, \bar y), \label{kt-yb} \nonumber\\
&&0= 
-r\widetilde\lambda \nabla_{2} f(\bar x, \bar u)
+\sum_{i\in I(\bar x,\bar u)} \mu_i^g \nabla_{2} g_i(\bar x, \bar u)+  [\nabla_{22}^2 f + \sum_{i\in I(\bar x,\bar u)} \bar{s}_i \nabla_{22}^2 g_i](\bar x, \bar u)^T\mu_{\phi},\label{kt-ub}\nonumber\\
&& 0= \nabla_{2} g_i(\bar x, \bar u)^T\mu_{\phi} -\mu_i^s,\quad i=1,\cdots,m,\label{kt-sb0}\\
&&\mu_i^g=0, i\in I_s(\bar x,\bar u),\
\mu_i^s=0,\ i\in I_g(\bar x,\bar u),\  
\mu_i^g>0,\ \mu_i^s>0,\ {\rm or}\ \mu_i^g\mu_i^s= 0,\ i\in I_0(\bar x,\bar u).\label{kt-sb}
\end{eqnarray}
The MPEC-NNAMCQ for $\widetilde{\cal{F}}$ imply that $r>0$ and we assume $r=1$ without loss of generality.

From the definitions of the index sets and the conditions (\ref{kt-sb0})-(\ref{kt-sb}), we have that
 \begin{eqnarray*}
&&\mu_i^g=0, i\in I_s(\bar x,\bar u),\
\nabla_{2} g_i(\bar x, \bar u)^T\mu_{\phi}=0,\ i\in I_g(\bar x,\bar u),\\  
&& \mu_i^g>0,\ \nabla_{2} g_i(\bar x, \bar u)^T\mu_{\phi}>0,\ {\rm or}\ \mu_i^g\nabla_{2} g_i(\bar x, \bar u)^T\mu_{\phi}= 0,\ i\in I_0(\bar x,\bar u).
\end{eqnarray*}
The conclusion holds.
\BOX

\subsection{Stationarity Relationships with Existing Formulations}
This subsection compares the stationarity conditions of SVF with three existing formulations: 
KKT-based  formulation (KP), Value function formulation (VP)  and combined problem (CP).
%

We first demonstrate that an S-stationary point of the problem SVF reduces to an S-stationary point of KP when $\bar y=\bar u$.
\begin{thm}\label{optmp}
Suppose that a feasible point $(\bar x,\bar y,\bar y,\bar s)$ is an S-stationary point of the problem $SVF$.
Then $(\bar x,\bar y,\bar s)$ is an S-stationary point of $KP$.
\end{thm}
{\bf Proof.} Since $(\bar x,\bar y,\bar y,\bar s)$ is feasible for $SVF$, then $(\bar x,\bar y, \bar s)$ is a feasible point of KP. 
If $(\bar x,\bar y,\bar y,\bar s)$ is an S-stationary point of the problem $SVF$,
there exist  $\lambda\in\mathbb{R}_+^{m+1}$, $\mu^{\phi}\in \mathbb{R}^l, \mu^g\in \mathbb{R}^m$  such that (\ref{lkt-xb})-(\ref{lkt-sb1}) and (\ref{lkt-sbs}) hold with $\bar y=\bar u$.
Set $\mu_i:=\mu_i^g+\lambda_i$, we have that
\begin{eqnarray}
&&0= \nabla_{1} F(\bar x, \bar y)  
+ [\nabla_{21}^2 f + \sum_{i\in I(\bar x,\bar y)} \bar{s}_i \nabla_{21}^2 g_i
](\bar x, \bar y)^T\mu_{\phi}
+ \sum_{i\in I(\bar x,\bar y)} \mu_i \nabla_{1} g_i(\bar x, \bar y), \label{mkt-xb}\nonumber \\
&&0=
\nabla_{2} F(\bar x, \bar y) 
+ \sum_{i\in I(\bar x,\bar y)} \mu_i \nabla_{2} g_i(\bar x, \bar y)
+  [\nabla_{22}^2 f + \sum_{i\in I(\bar x,\bar y)} \bar{s}_i \nabla_{22}^2 g_i](\bar x, \bar y)^T\mu_{\phi},\label{mkt-yb}\\
&&\mu_i=0, i\in I_s(\bar x,\bar y),\quad
\nabla_{2} g_i(\bar x, \bar y)^T\mu_{\phi}=0,\ i\in I_g(\bar x,\bar y),\label{mkt-sb1}\nonumber\\
&&\mu_i\geq 0,\ \nabla_{2} g_i(\bar x, \bar y)^T\mu_{\phi}\geq 0,\ i\in I_0(\bar x,\bar y),\label{mkt-sb2}\nonumber
\end{eqnarray}
where (\ref{mkt-yb}) is obtained by summing the conditions (\ref{lkt-yb}) and (\ref{lkt-ub}).
Then $(\bar x,\bar y, \bar s)$ is an S-stationary point of KP. 
\BOX

To construct the relationship of stationary points among the problems SVF and VP, CP, we assume $(\bar u,\bar s)$ satisfies 
\begin{eqnarray}\label{subv}
\nabla_{1} f(\bar x, \bar u)+\sum_{i\in I(\bar x,\bar u)} \bar{s}_i \nabla_{1} g_i(\bar x, \bar u)\in \partial^c V(\bar x).
\end{eqnarray} 
Let  ${\cal{F}}(\bar x):=\{(u,s):\nabla_{2} f(\bar x,u)+s^T \nabla_{2} g(\bar x,u)=0,\ 0\geq g(\bar x,u)\perp s\geq 0\}$. 

\begin{thm}\label{optvp}
(i) Assume the Assumption \ref{pseudoconvex} holds and $(\bar x,\bar y,\bar u,\bar s)\in {\cal{F}}_S$ is an M-stationary point of $SVF$ such that (\ref{subv}) holds.
If the MPEC-NNAMCQ for the system ${\cal{F}}(\bar x)$ holds at $(\bar u,\bar s)$, then $(\bar x,\bar y)$ is a stationary point of $VP$.\\
(ii) Suppose that $(\bar x,\bar y)$ is a stationary point of VP.
If $f(x,y)$ is jointly convex, $g(x,y)$ is jointly quasiconvex,
and $M(\bar x,\bar u)=\{\bar s\}$,
then $(\bar x,\bar y,\bar u,\bar s)$ is an S-stationary point of SVF.
\end{thm}
{\bf Proof.} (i) From the Assumption \ref{pseudoconvex}, $(\bar x,\bar y)$ is feasible for VP.
Since $(\bar x,\bar y,\bar u,\bar s)$ is an M-stationary point of the problem $SVF$, there exist  $\lambda\in\mathbb{R}_+^{m+1}$, $\mu^{\phi}\in \mathbb{R}^l, \mu^g\in \mathbb{R}^m$  such that $\lambda_0(f(\bar x, \bar y)-f(\bar x, \bar u))=0$ and (\ref{lkt-xb})-(\ref{lkt-sb1}), (\ref{lkt-sbm}) hold.

From the feasibility of $(\bar x,\bar y,\bar u,\bar s)$ for SVF,  $-\nabla_{2} f(\bar x, \bar u)=\sum_{i\in I(\bar x,\bar u)} \bar{s}_i \nabla_{2} g_i(\bar x, \bar u)$ and thus the conditions (\ref{lkt-xb}) and (\ref{lkt-ub}) reduce to
\begin{eqnarray}
&&0=  \nabla_{1} F(\bar x, \bar y)+
\lambda_0(\nabla_{1} f(\bar x, \bar y)-\nabla_{1} f(\bar x, \bar u)-\sum_{i\in I(\bar x,\bar u)} \bar{s}_i \nabla_{1} g_i(\bar x, \bar u) )
+\sum_{i\in I(\bar x,\bar y)} \lambda_i \nabla_{1} g_i(\bar x, \bar y) \nonumber\\
&&~~~+ [\nabla_{21}^2 f + \sum_{i\in I(\bar x,\bar u)} \bar{s}_i \nabla_{21}^2 g_i
](\bar x, \bar u)^T\mu_{\phi}
+ \sum_{i\in I(\bar x,\bar u)} (\mu_i^g+\lambda_0 \bar{s}_i) \nabla_{1} g_i(\bar x, \bar u),
 \label{kt-xb1} \\
&&0= 
\sum_{i\in I(\bar x,\bar u)} (\mu_i^g+\lambda_0 \bar{s}_i) \nabla_{2} g_i(\bar x, \bar u)+  [\nabla_{22}^2 f + \sum_{i\in I(\bar x,\bar u)} \bar{s}_i \nabla_{22}^2 g_i](\bar x, \bar u)^T\mu_{\phi}.\label{kt-ub1}
\end{eqnarray}
Denote by $\tilde{\mu}_i^g:=\mu_i^g+\lambda_0 \bar{s}_i$ for each $i\in I(\bar x,\bar u)$.
For $ i\in I_0(\bar x,\bar u)$, we have that $\bar{s}_i=0$. Followed from (\ref{lkt-sbm}), we have that
\begin{eqnarray}\label{lkt-sbm-1}
\tilde{\mu}_i^g> 0,\ \nabla_2 g_i(\bar x, \bar u)^T\mu_{\phi}> 0,\ {\rm or}\ \tilde{\mu}_i^g \nabla_2 g_i(\bar x, \bar u)^T\mu_{\phi}=0,\ i\in I_0(\bar x,\bar u).
\end{eqnarray}

From MPEC-NNAMCQ for the system ${\cal{F}}(\bar x)$, the conditions (\ref{lkt-sb1}), (\ref{kt-ub1}) and (\ref{lkt-sbm-1})
imply that $\mu_{\phi}= 0$, $\tilde{\mu}_i^g=0$ for each $i$.
Followed from (\ref{subv}), the condition (\ref{kt-xb1}) reduces to
\begin{eqnarray}
0&=&\nabla_{1} F(\bar x, \bar y)+\lambda_0(\nabla_{1} f(\bar x,\bar y)-\nabla_{1} f(\bar x,\bar u)-\sum_{i\in I(\bar x,\bar u)} \bar{s}_i \nabla_{1} g_i
(\bar x,\bar u))+\sum_{i\in I(\bar x,\bar y)} \lambda_i \nabla_{1} g_i(\bar x,\bar y)\nonumber\\
&\in &\nabla_{1} F(\bar x, \bar y)+\lambda_0(\nabla_{1} f(\bar x,\bar y)-\partial^c V(\bar x)) 
+\sum_{i\in I(\bar x,\bar y)} \lambda_i \nabla_{1} g_i(\bar x,\bar y).\label{nvp}
\end{eqnarray}
Then the conditions (\ref{lkt-yb}) and (\ref{nvp}) imply that 
 $(\bar x,\bar y)$ is a stationary point of $VP$.
 
 (ii) 
It is easy to see $(\bar x,\bar y,\bar u,\bar s)\in {\cal{F}}_S$. 
 From the Proposition \ref{partialV} (iii), $\partial^c V(\bar x)=\nabla_1 f(\bar x,\bar u)+\bar{s}^T \nabla_1 g(\bar x,\bar u)$.
Then there exist $\mu\geq 0$, $\lambda\in\mathbb{R}_+^m$ such that 
\begin{eqnarray*}\label{stavp}
0= \nabla F(\bar x,\bar y) +\mu [\nabla f(\bar x,\bar y) -(\nabla_1 f(\bar x,\bar u)+\bar{s}^T \nabla_1 g(\bar x,\bar u))\times \{0\}_l] +\sum_{I(\bar x,\bar y)} \lambda_i \nabla g_i(\bar x,\bar y).
\end{eqnarray*}
By setting $\mu_{\phi}:= 0$, $\mu_i^g:=-\mu\bar {s}_i$ for $i=1,\cdots,m$, we observe that  $(\bar x,\bar y,\bar u,\bar s)$ is an S-stationary point of the problem $SVF$. We complete the proof.
 \BOX

We now recall the combined problem \cite{yz2}:
\begin{eqnarray*}
({\rm CP})~~~~~~\min_{x,y,s} && F(x,y)\nonumber\\
{\rm s.t.} && f(x,y)-V(x) \leq 0,\\
&&  \nabla_2 f(x,y)+ \sum_{i=1}^m s_i \nabla_{2} g_i(x,y) =0, \label{FC}\\
&& s_i\geq 0,\ g_i(x,y)\leq 0,\ s_i g_i(x,y)=0,\ i=1,\cdots,m.\nonumber
\end{eqnarray*}
In the rest of this section, we investigate  relationship between stationary points of SVF and CP.
\begin{thm}\label{optcp}
Assume the Assumption \ref{pseudoconvex} holds
 and $(\bar x,\bar y,\bar y,\bar s)\in {\cal{F}}_S$ is an S-stationary point of $SVF$.
 If  (\ref{subv}) holds with $\bar y=\bar u$, then $(\bar x,\bar y,\bar s)$ is an S-stationary point of $CP$.
\end{thm}
{\bf Proof.} From the Assumption \ref{pseudoconvex}, $(\bar x,\bar y,\bar s)$ is a feasible point of $CP$.
Since $(\bar x,\bar y,\bar y,\bar s)\in {\cal{F}}_S$ is an S-stationary point of $SVF$, 
there exist  $\lambda\in\mathbb{R}_+^{m+1}$, $\mu^{\phi}\in \mathbb{R}^l, \mu^g\in \mathbb{R}^m$  such that (\ref{lkt-xb})-(\ref{lkt-sb1}) and (\ref{lkt-sbs}) hold with $\bar y=\bar u$.

Similarly as Theorem \ref{optvp}, (\ref{lkt-xb}) and (\ref{lkt-ub}) reduce to (\ref{kt-xb1}) and (\ref{kt-ub1}) with $\bar y=\bar u$.
Denote by $\nu_i^g:=\mu_i^g+\lambda_0 \bar{s}_i+\lambda_i$ for each $i\in I(\bar x,\bar u)$.
Summing the conditions (\ref{lkt-yb}) and (\ref{kt-ub1}), we have that
\begin{eqnarray*}0= 
\nabla_2 F(\bar x, \bar y)+\lambda_0\nabla_2 f(\bar x, \bar y)+\sum_{i\in I(\bar x,\bar y)} \nu_i^g \nabla_{2} g_i(\bar x, \bar y)+  [\nabla_{22}^2 f + \sum_{i\in I(\bar x,\bar y)} \bar{s}_i \nabla_{22}^2 g_i](\bar x, \bar y)^T\mu_{\phi}.\label{kt-ub2}
\end{eqnarray*}
Followed from (\ref{subv}) with $\bar y=\bar u$, the condition (\ref{kt-xb1}) reduces to
\begin{eqnarray*}
0&\in &\nabla_{1} F(\bar x, \bar y)+\lambda_0(\nabla_{1} f(\bar x,\bar y)-\partial^c V(\bar x))
+ \sum_{i\in I(\bar x,\bar y)} \nu_i^g \nabla_{1} g_i(\bar x, \bar y)\nonumber \\
&&~~+ [\nabla_{21}^2 f + \sum_{i\in I(\bar x,\bar y)} \bar{s}_i \nabla_{21}^2 g_i
](\bar x, \bar y)^T\mu_{\phi}.\label{nvp-1}
\end{eqnarray*}
Similarly as the proof of Theorem \ref{optmp},
 $(\bar x,\bar y,\bar s)$ is an S-stationary point of $CP$ from \cite[Definition 3.2]{yz2}.
\BOX

Note that the condition (\ref{subv}) holds automatically when the conditions in Proposition \ref{partialV} (iii) hold.
Then the condition (\ref{subv}) in Theorems \ref{optvp} (i), \ref{optcp} can be replaced by the jointly convexity of $f(x,y)$ and  jointly quasiconvexity of $g(x,y)$.
 We summarize the relationships among various stationary points  in Figure 2.
\begin{figure}%
\centering
	\scriptsize
		 \tikzstyle{format}=[rectangle,draw,thin,fill=white]
		
		 \tikzstyle{test}=[diamond,aspect=2,draw,thin]
		
		\tikzstyle{point}=[coordinate,on grid,]
\begin{tikzpicture}
        \node[format](SVF C-stationary){SVF C-stationary};
        \node[format,right of=SVF C-stationary,node distance=32mm](SVF M-stationary){SVF M-stationary};
        \node[format,right of=SVF M-stationary,node distance=32mm](SVF S-stationary){SVF S-stationary};
       \node[format,right of=SVF S-stationary,node distance=32mm](KP S-stationary){KP S-stationary};
        \node[format, below of=SVF S-stationary,node distance=13mm](CP S-stationary){CP S-stationary};
        \node[format, above of=SVF M-stationary,node distance=13mm](VP stationary){VP stationary};

\draw[->](SVF M-stationary)--(SVF C-stationary);
\draw[->](SVF M-stationary)--(VP stationary);
\draw[->](SVF S-stationary)--(SVF M-stationary);
\draw[->](SVF S-stationary)--(CP S-stationary);
\draw[->](SVF S-stationary)--(KP S-stationary);
\draw[->](VP stationary)--(SVF S-stationary);
\draw[->] (VP stationary) -- node[right] {jointly convex + LICQ}(SVF S-stationary);
\draw[->] (SVF S-stationary) -- node[above] {$\bar y=\bar u$}(KP S-stationary);
\draw[->] (SVF S-stationary) -- node[left] {$\bar y=\bar u$ }(CP S-stationary);
\draw[->] (SVF S-stationary) -- node[right] {(\ref{subv})}(CP S-stationary);
\draw[->] (SVF M-stationary) -- node[right] {(\ref{subv})}(VP stationary);
\draw[->] (SVF M-stationary) -- node[left] {CQ}(VP stationary);
\end{tikzpicture}

\centering{Fig.2 Relation among stationary points}
\end{figure}

\section{Smoothing‑Barrier Augmented Lagrangian Approximation}
To handle the complementarity constraints, the SVF model is approximated by a family of smooth subproblems using the SBAL method \cite{ld}. It is shown that, as the smoothing parameter tends to zero, the solution sets of the smooth problems approach that of the SVF model, and their stationary points converge to the C-stationary points of the SVF formulation, provided that the Lagrange multipliers remain uniformly bounded.


For $r>0$, by introducing auxiliary variables $z_i> 0, i=1,\cdots,m$, $({\rm P}_x)$ is approximated by 
 \begin{eqnarray*}
({\rm P}_x^{r})~~~~~~~~\min_{u,z} &&  f(x,u)-r \sum_{i=1}^m \ln z_i\\
{\rm s.t.} && z_i+g_i(x,u)=0, i=1,\cdots,m.
\end{eqnarray*}
For each $\rho>0$, $s\in\mathbb{R}^m$, 
the augmented Lagrangian function of $({\rm P}_x^{r})$ is defined as follows:
\begin{eqnarray*}
\varphi_r^{\rho}(x,u,z,s):= f(x,u)+ \sum_{i=1}^m [-r \ln z_i+ s_i(z_i+g_i(x,u))+\frac{1}{2 \rho}(z_i+g_i(x,u))^2].
\end{eqnarray*}
To ensure that $s$ is a good estimate of Lagrange multiplier vector, the authors in \cite{ld} maximized the augmented Lagrangian with respect to $s$. This leads to the following problem:
\begin{eqnarray*}\label{minmax}
({\rm P}_x^{r,\rho})~~~~~~~~
\min_{u,z} \max_s  &&  \varphi_r^{\rho}(x,u,z,s).
\end{eqnarray*}

For each $x$ and any solution $(u,z,s)$ of the problem $({\rm P}_x^{r,\rho})$, we must have that $\nabla_z \varphi_r^{\rho}(x,u,z,s)=0$, 
which derives that $z_i$ is a function depending on $x,u,s_i,r,\rho$,
 \begin{eqnarray}\label{z}
z_i(x,u,s_i,r,\rho):= \frac{1}{2}[\sqrt{(\rho s_i+ g_i(x,u))^2+4 r \rho}-(\rho s_i+g_i(x,u))],\ i=1,\cdots,m.
\end{eqnarray}
Since the function $\varphi_r^{\rho}$ is convex with respect to the variable $z$, we can replace $z$ by $z(x,u,s,r,\rho)$ 
and rewrite $\varphi_r^{\rho}(x,u,z,s)$ by $f_r^{\rho}(x,u,s)$. 
For any $x$ and $r>0, \rho>0$,
any solution of $({\rm P}_x^{r,\rho})$ satisfies the KKT conditions, 
\begin{eqnarray}
&&\phi^{r,\rho}(x,u,s):=\nabla_{2} f_r^{\rho}(x,u,s)=\nabla_{2} f(x,u)+  \sum_{i=1}^m \frac{\kappa_i(x,u,s_i,r,\rho)}{\rho} \nabla_{2} g_i(x,u) =0,\label{kkt1}\\
&&\psi^{r,\rho}(x,u,s):=\nabla_s f_r^{\rho}(x,u,s)=z(x,u,s,r,\rho)+g(x,u)=0,\label{kkt2}
\end{eqnarray}
where for $i=1,\cdots,m$, 
\begin{eqnarray}\label{k}
\kappa_i(x,u,s_i,r,\rho)&:=&z_i(x,u,s_i,r,\rho)+g_i(x,u) +\rho s_i\nonumber\\
&=&\frac{1}{2}[\sqrt{(\rho s_i+ g_i(x,u))^2+4 r \rho}+(\rho s_i+g_i(x,u))].
\end{eqnarray}
We write $z_i(x,u,s_i,r,\rho)$ and $\kappa_i(x,u,s_i,r,\rho)$ as $z_i$ and $\kappa_i$ for convince, respectively.
Similarly as Liu et al. \cite{ld}, properties for $z_i$ and $\kappa_i$ are stated in the following lemma.

\begin{lemma}\label{zkproperty}
For each $r\geq 0, \rho> 0$,  $i=1,\cdots,m$, the following conclusions hold.\\
(1)  $z_i\geq 0$, $\kappa_i\geq 0$, $z_i+g_i(x,u)=\kappa_i-\rho s_i$ and $ z_i \kappa_i=\rho r$;\\
(2)   $g_i(x,u)+z_i=0$ if and only if $g_i(x,u)\leq 0$, $s_i\geq 0$ and $s_i g_i(x,u)=-r$.\\
(3) For $r>0$, $z_i$ and $\kappa_i$ are differentiable with respect to the variable $x, u$, $s$ and $\rho$, respectively,
\begin{eqnarray*}
 \nabla_{12} z_i=\frac{- z_i}{z_i+\kappa_i} \nabla g_i(x,u),
 && \nabla_{12} \kappa_i=\frac{ \kappa_i}{z_i+\kappa_i} \nabla g_i(x,u),\\
 \nabla_s z_i=\frac{- \rho z_i}{z_i+\kappa_i}e_i,
 &&\nabla_s \kappa_i=\frac{\rho\kappa_i}{z_i+\kappa_i}e_i.
\end{eqnarray*}
\end{lemma}

By Lemma \ref{zkproperty}(1)–(2), for any penalty  parameter $\rho>0$, the perturbed KKT system (\ref{kkt1})-(\ref{kkt2}) provides a smooth approximation to the complementarity conditions (\ref{fl1})-(\ref{fl2}).
Consequently, SVF can be approximated by the following smoothed subproblems:
\begin{eqnarray*}
({\rm SVF}_{r,\rho})~~~~~~~\min_{x,y,u,s}  && F(x,y)\nonumber\\
{\rm s.t.} && f(x,y)-f(x,u)\leq 0,\label{ff0}\\
&&\phi^{r,\rho}(x,u,s)=\nabla_{2} f(x,u)+  \sum_{i=1}^m \frac{\kappa_i}{\rho} \nabla_{2} g_i(x,u) =0,\nonumber\\
&& \psi^{r,\rho}(x,u,s)=z(x,u,s,r,\rho)+g(x,u)=0,\nonumber\\
&& g(x,y)\leq 0.\nonumber
\end{eqnarray*}
The following theorem shows that as the smoothing parameter $r$ rends to zero,  the solution sets of the smoothed problems $({\rm SVF}_{r,\rho})$ converge to the solution set of  SVF for a fixed penalty parameter $\rho>0$.

Denote the feasible region and solution set of $({\rm SVF}_{r,\rho})$ by ${\cal{F}}_{r,\rho}$ and $S_{r,\rho}$, respectively.  
From \cite{ld}, the following assumption holds if for any $x$, there exists an interior point $u$, i.e., $g(x,u)<0$.
\begin{ass}\label{nempty}
For the point of interest  $x$ and any $r>0$, $\rho>0$, suppose $S_r^{\rho}(x)$ is nonempty, where  $S_r^{\rho}(x)$ denotes  the solution set of the system (\ref{kkt1})-(\ref{kkt2}).
\end{ass}

\begin{thm}\label{feavp}
Assume the Assumptions \ref{pseudoconvex}, \ref{nempty} hold.
Then for any $\bar{\rho}>0$,
\begin{eqnarray*}
\limsup_{r\to 0, \rho\to\bar{\rho}}S_{r,\rho} \subseteq S_{SV},
\end{eqnarray*} 
where $S_{SV}$ is the solution set of the problem SVF.
\end{thm}
 {\bf Proof.}   
 Let $(x_r^{\rho}, y_r^{\rho}, u_r^{\rho},s_r^{\rho})\in S_{r,\rho}$.
 Assume without loss of generality that $\lim_{r\to 0,\rho\to \bar{\rho}}(x_r^{\rho}, y_r^{\rho}, u_r^{\rho},s_r^{\rho})=(\bar x,\bar y,\bar u,\bar s)$. 
 Taking limits as $r\to 0,\rho\to \bar{\rho}$, from Lemma \ref{zkproperty} (1)-(2), it follows that
\begin{eqnarray}
&&\nabla_{2} f(\bar x,\bar u)+  \sum_{i=1}^m \bar{s}_i \nabla_{2} g_i(\bar x,\bar u) =0,\label{lkt1}\\
&&g_i(\bar x,\bar u)\leq 0,\ \bar{s}_i\geq 0,\ \bar{s}_i g_i(\bar x,\bar u)=0,\ i=1,\cdots,m,
\label{lkt2}\\
&&f(\bar x,\bar y)-f(\bar x,\bar u)\leq 0,\quad g(\bar x,\bar y)\leq 0.\label{lkt3}
\end{eqnarray}
Then $(\bar x,\bar y,\bar u,\bar s)\in \cal{F}_S$.


For any $(x^*,y^*,u^*,s^*)\in \cal{F}_S$,  similar to (\ref{vl}), we have that $y^*\in S(x^*)$.
From the  Assumption \ref{nempty}, there exists $(\tilde{u}_r^{\rho},\tilde{s}_r^{\rho})\in S_r^{\rho}(x^*)$ and thus  $(x^*,y^*,\tilde{u}_r^{\rho},\tilde{s}_r^{\rho})\in {\cal{F}}_{r,\rho}$. Then we have that $F(x^*,y^*)\geq F(x_r^{\rho}, y_r^{\rho})$ since $(x_r^{\rho}, y_r^{\rho}, u_r^{\rho},s_r^{\rho})\in S_{r,\rho}$.
Taking limits as $r\to 0,\rho\to \bar{\rho}$, $F(x^*,y^*)\geq F(\bar x,\bar y)$ and thus
 $(\bar x,\bar y,\bar u,\bar s)$ is a solution of SVF.
\BOX

In the rest of this section, we prove that any accumulation point of the stationary points of $({\rm SVF}_{r,\rho})$ is a C-stationary point of SVF.
\begin{thm}\label{stationary} 
Suppose that $(x_r^{\rho}, y_r^{\rho}, u_r^{\rho},s_r^{\rho})$ is a stationary point of $({\rm SVF}_{r,\rho})$ for each $r>0, \rho>0$. 
Assume that the corresponding multipliers $\{(\lambda^{r,\rho},\mu^{r,\rho},\beta^{r,\rho})\}$ are uniformly bounded,
then any accumulation point of $\{(x_r^{\rho}, y_r^{\rho}, u_r^{\rho},s_r^{\rho})\}$, as $r\to 0, \rho\to\bar{\rho}$, $\bar{\rho}>0$, is a C-stationary point of SVF.
\end{thm}
{\bf Proof.}  Let $z^{r,\rho}:=z(x_r^{\rho},  u_r^{\rho},s_r^{\rho},r,\rho)$, $\kappa^{r,\rho}:=\kappa(x_r^{\rho},  u_r^{\rho},s_r^{\rho},r,\rho)$. 
Without loss of generality, we suppose there exist $\bar{\kappa},\bar{z}$, $(\bar x,\bar y,\bar{u},\bar s)$ such that 
$\bar{\kappa}=\displaystyle\lim_{r\to 0, \rho\to\bar{\rho}} \kappa^{r,\rho}$, $\bar{z}=\displaystyle\lim_{r\to 0, \rho\to\bar{\rho}} z^{r,\rho}$, $\displaystyle\lim_{r\to 0, \rho\to\bar{\rho}} (x_r^{\rho}, y_r^{\rho}, u_r^{\rho},s_r^{\rho})=(\bar x,\bar y,\bar{u},\bar s)$.
Similarly with the proof of Theorem \ref{feavp}, $(\bar x,\bar y,\bar{u}, \bar{s})$ is a feasible point of SVF.

Since $\{(x_r^{\rho}, y_r^{\rho}, u_r^{\rho},s_r^{\rho})\}$ is a stationary point for the problem $({\rm SVF}_{r,\rho})$,
\begin{eqnarray}
0&=& \nabla_{1} F(x_r^{\rho},y_r^{\rho})  
+[\nabla_{12}^2 f + \sum_{i=1}^m \frac{\kappa_i^{r,\rho}}{\rho} \nabla_{12}^2 g_i](x_r^{\rho},u_r^{\rho})^T \mu_r^{\rho}
\nonumber \\
&+&\sum_{i=1}^m \frac{\kappa_i^{r,\rho} \nabla_{1} g_i(x_r^{\rho},u_r^{\rho})}{\rho(z_i^{r,\rho}+\kappa_i^{r,\rho})}  \nabla_{2} g_i(x_r^{\rho},u_r^{\rho})^T \mu_r^{\rho}
+\sum_{i=1}^m \frac{\beta_i^{r,\rho} \kappa_i^{r,\rho}}{z_i^{r,\rho}+\kappa_i^{r,\rho}}\nabla_{1} g_i(x_r^{\rho},u_r^{\rho})
\nonumber\\
&+& \lambda_0^{r,\rho}(\nabla_{1} f(x_r^{\rho},y_r^{\rho})-\nabla_{1} f(x_r^{\rho},u_r^{\rho}))
+\sum_{i=1}^m \lambda_i^{r,\rho} \nabla_{1} g_i(x_r^{\rho},y_r^{\rho}),
\label{alfp0kp-1al}\\
0&=&\nabla_{2} F(x_r^{\rho},y_r^{\rho}) +\lambda_0^{r,\rho} \nabla_{2} f(x_r^{\rho},y_r^{\rho})+\sum_{i=1}^m \lambda_i^{r,\rho} \nabla_{2} g_i(x_r^{\rho},y_r^{\rho}), \label{alfp1kp-al}\\
0&=&-\lambda_0^{r,\rho}\nabla_{2} f(x_r^{\rho},u_r^{\rho})+\sum_{i=1}^m \frac{\kappa_i^{r,\rho}}{\rho(z_i^{r,\rho}+\kappa_i^{r,\rho})} \nabla_{2} g_i(x_r^{\rho},u_r^{\rho}) \nabla_{2} g_i(x_r^{\rho},u_r^{\rho})^T \mu_r^{\rho}\nonumber\\
&+&  \sum_{i=1}^m \frac{\beta_i^{r,\rho} \kappa_i^{r,\rho}}{z_i^{r,\rho}+\kappa_i^{r,\rho}}\nabla_{2} g_i(x_r^{\rho},u_r^{\rho})+ [\nabla_{22}^2 f + \sum_{i=1}^m \frac{\kappa_i^{r,\rho}}{\rho} \nabla_{22}^2 g_i](x_r^{\rho},u_r^{\rho})^T \mu_r^{\rho},
\label{alfp11kp-1al}
\\0&=& \frac{\kappa_i^{r,\rho} \nabla_{2} g_i(x_r^{\rho},u_r^{\rho})^T \mu_r^{\rho}}{z_i^{r,\rho}+\kappa_i^{r,\rho}} 
-\frac{\rho \beta_i^{r,\rho} z_i^{r,\rho}}{z_i^{r,\rho}+\kappa_i^{r,\rho}}.\label{alfp1vkp-1al}
\end{eqnarray}
From (\ref{alfp1vkp-1al}), the conditions (\ref{alfp0kp-1al}), (\ref{alfp11kp-1al}) reduce to
\begin{eqnarray}
0&=& \nabla_{1} F(x_r^{\rho},y_r^{\rho})  
+[\nabla_{12}^2 f + \sum_{i=1}^m \frac{\kappa_i^{r,\rho}}{\rho} \nabla_{12}^2 g_i](x_r^{\rho},u_r^{\rho})^T \mu_r^{\rho}
+\sum_{i=1}^m \beta_i^{r,\rho}\nabla_{1} g_i(x_r^{\rho},u_r^{\rho})
\nonumber\\
&+& \lambda_0^{r,\rho}(\nabla_{1} f(x_r^{\rho},y_r^{\rho})-\nabla_{1} f(x_r^{\rho},u_r^{\rho}))
+\sum_{i=1}^m \lambda_i^{r,\rho} \nabla_{1} g_i(x_r^{\rho},y_r^{\rho}),
\label{alfp0kp-1al-1}\\
0&=&-\lambda_0^{r,\rho}\nabla_{2} f(x_r^{\rho},u_r^{\rho})
+\sum_{i=1}^m \beta_i^{r,\rho}\nabla_{2} g_i(x_r^{\rho},u_r^{\rho})
 +  [\nabla_{22}^2 f + \sum_{i=1}^m \frac{\kappa_i^{r,\rho}}{\rho} \nabla_{22}^2 g_i](x_r^{\rho},u_r^{\rho})^T \mu_r^{\rho}.
\label{alfp11kp-1al-1}
\end{eqnarray}

We assume there exists $(\lambda_0,\lambda,\mu,\beta)$ such that $\displaystyle\lim_{r\to 0, \rho\to\bar{\rho}} (\lambda_0^{r,\rho}, \lambda^{r,\rho},\mu^{r,\rho},\beta^{r,\rho})=(\lambda_0,\lambda,\mu,\beta)$.
If $i\notin I(\bar x,\bar y)$, from the continuity and complementary, $g_i(x_r^{\rho},y_r^{\rho})<0$ and $ \lambda_i^{r,\rho}=0$, thus $\lambda_i= 0$. 
If $i\notin I(\bar x,\bar u)$, we have that  $\bar{z}_i>0$ and thus $\bar{\kappa}_i=0$ from the Lemma \ref{zkproperty} (1) and the fact that $r\to 0$.  Moreover, $\beta_i=0$ from (\ref{alfp1vkp-1al}) and the continuity in this case.

Taking limits on both sides of (\ref{alfp1kp-al}), (\ref{alfp0kp-1al-1})-(\ref{alfp11kp-1al-1}) as $r\to 0, \rho\to\bar{\rho}$, we have that
\begin{eqnarray}
0&=& 
\nabla_{1} F(\bar x,\bar y)  +\lambda_0(\nabla_{1} f(\bar x,\bar y)-\nabla_{1} f(\bar x,\bar u))
+\sum_{i\in I(\bar x,\bar y)} \lambda_i \nabla_{1} g_i(\bar x,\bar y)\nonumber\\
&&+[\nabla_{12}^2 f + \sum_{i\in I(\bar x, \bar u)} \frac{\bar{\kappa}_i}{\bar{\rho}} \nabla_{12}^2 g_i](\bar x,\bar u)^T \mu
+
\sum_{i\in I(\bar x,\bar u)}\beta_i \nabla_{1} g_i(\bar x,\bar u),
\label{alktx1} \\
0&=&\nabla_{2} F(\bar x,\bar y) +\lambda_0 \nabla_{2} f(\bar x,\bar y)+\sum_{i\in I(\bar x,\bar y)} \lambda_i \nabla_{2} g_i(\bar x,\bar y), \label{alkty}\\
0&=& 
-\lambda_0 \nabla_{2} f(\bar x,\bar u )+
\sum_{i\in I(\bar x,\bar u)} \beta_i \nabla_{2} g_i(\bar x,\bar u)+  [\nabla_{22}^2 f + \sum_{i\in I(\bar x, \bar u)} \frac{\bar{\kappa}_i}{\bar{\rho}} \nabla_{22}^2 g_i](\bar x,\bar u)^T \mu. \label{alktu1}
\end{eqnarray}
For $i\in I_g(\bar x,\bar u)$, we have that $\bar{\kappa}_i>0$ and $\bar{z}_i=0$.
For $i\in I_s(\bar x,\bar u)$, we have that $\bar{\kappa}_i=0$ and $\bar{z}_i>0$.
Taking limits on both sides of (\ref{alfp1vkp-1al}) for such $i$, we have that
\begin{eqnarray}\label{alkts1}
\beta_i=0,\ i\in I_s(\bar x,\bar u),\quad 
\nabla_{2} g_i(\bar x,\bar u)^T \mu=0,\ i\in I_g(\bar x,\bar u).
\end{eqnarray}
For $i\in I_0(\bar x,\bar u)$, combining condition (\ref{alfp1vkp-1al}) with $z^{r,\rho} \ge 0$ and $\kappa^{r,\rho} \ge 0$ implies that
\begin{eqnarray}\label{l3}
	\beta_i \nabla_{2} g_i(\bar x, \bar u)^T \mu \geq 0,i\in I_0(\bar x,\bar u).
\end{eqnarray}

The conditions (\ref{alktx1})-(\ref{l3}) imply that $(\bar x,\bar y,\bar{u},\bar s)$ is a C-stationary point of SVF.						
We complete the proof.
\BOX

\section{Numerical experiments}
This section evaluates the performance of the proposed SVF model. We compare the SVF solved via SBAL with three representative methods: the KKT model with a relaxation method, the SBAL-based smoothing algorithm applied to KP
and the extension of Mond-Weir duality reformulation with a relaxation method (eRMDP) \cite{llz}.
The tests covered 109 nonlinear bilevel problems from the BOLIB library \cite{bolib} with constrained lower level problems.  All experiments were implemented in MATLAB\_R2021a and run on a MacBook Pro  (Apple M1 chip, 512GB SSD).


\subsection{Algorithmic Framework}
We first summarize the implementation of the SBAL-based algorithm for solving the SVF model (SVF-SBAL). The algorithm proceeds by solving a sequence of smooth approximations to the original complementarity constrained problem, which is detailed below.
\begin{alg}\label{algo3-1}
\begin{enumerate} 
\item Given  initial point $(x_0,y_0,u_0,s_0)\in\mathbb{R}^{d+2l+m}$,  
initial parameters $r_0>0$, $\rho_0>0$.
Let $\bar{\rho}>0$ and \( \delta\) be a constant in $[0,1)$.
Set \(k:=0\).


\item[{2.}] If a stopping criterion leading to the stationary conditions for SVF hold, terminate. 

\item[{3.}] Find an approximated stationary point $\{(x_{k+1},y_{k+1},u_{k+1},s_{k+1})\}$ of the problem $({\rm SVF}_{r_k,\rho_k})$.
										
\item[{4.}]   Set $r_{k+1}:=\delta r_k$, $\rho_{k+1}:= \max\{\bar{\rho}, \delta \rho_k\}$ and go to Step 2.
\end{enumerate}
\end{alg}

As shown in Theorem \ref{nna-fail}, the MPEC-NNAMCQ condition fails at all feasible points of SVF, meaning the boundedness of multipliers may not be guaranteed. One may instead invoke weaker conditions such as MPEC relaxed constant positive linear dependence constraint qualification (MPEC-RCPLD) to establish convergence of stationary points without bounded multipliers; see \cite{xy2}. For clarity and focus, we omit such technical extensions in this paper.

For comparison, we also apply the same SBAL-based smoothing framework to KP; named as KP-SBAL.  
For the KKT relaxed method (KP-RLX), a relaxation technique is employed \cite{llz}, in which the complementarity constraints are relaxed by introducing a small positive scalar $\varepsilon$, replacing $s^{\top} g(x,y) = 0$ with $s^{\top} g(x,y) \leq \varepsilon$.  Similarly, in the eRMDP approach, the condition $f(x,y) \leq f(x,u)$ is relaxed to $f(x,y) - f(x,u) \leq \varepsilon$ \cite{llz}.
	In both cases, as $\varepsilon \to 0$, the feasible set of the relaxed problem approaches that of the original reformulation. 

\subsection{Experimental Setup}					
In Algorithm \ref{algo3-1}, the parameters were set as follows:
$\rho_0=1$, $r_0=1$, $\bar{\rho}=0.01$, $\delta=0.1$. The same parameter settings were applied to KP-SBAL. For other compared methods KP-RLX, eRMDP, we set the initial relaxation parameter to $\varepsilon_0 = 1$ and update it after calculation of each subproblem as $\varepsilon_{k+1}=0.1\varepsilon_k$.
The initial point $(x_0, y_0)$ was taken from the reference points provided by the BOLIB library and $s_0 \sim \mathcal{U}(0,1)$.
For the SVF-SBAL and eRMDP approaches, we set $u_0: = y_0$ for convex lower level problems. For nonconvex lower level problems, complementary strategies of $u_0: = y_0$ and  $u_0 := -y_0$ were tested. The algorithm was considered successful if either strategy converged, and the corresponding result was adopted.
%

The stopping criteria were defined as follows:
\begin{enumerate}
\item $\text{Res}_k\le5\times10^{-5}$;
\item  the iteration count $k\ge50$;
\item  the iteration count $k\ge20$ and $\|\text{Res}_k -\text{Res}_{k-1}\|\le 10^{-8}$;
\item  the iteration count $k\ge30$ and $\text{Res}_k\le5\times10^{-4}$.
\end{enumerate}
Here, $\text{Res}_k$ denotes the KKT residual norm of each reformulated model at iteration $k$.

Let $\mathcal{SL} = \{\text{SVF-SBAL}, \text{KP-SBAL}, \text{KP-RLX}, \text{eRMDP}\}$ be the solver set and $\mathcal{P}$ be the test problems set with $n_p = 109$. 
For each problem $p \in \mathcal{P}$ and solver $sl \in \mathcal{SL}$, we define the performance metric of solution error $\epsilon^x_{p,sl}$,  objective error $\epsilon^f_{p,sl}$ and composite error $\omega_{p,sl}$ by
$$
\epsilon^x_{p,sl}=\|(x_{p,sl}, y_{p,sl}) - (x_p^*,y_p^*)\|,\quad \epsilon^f_{p,sl}=\|((F_{p,sl}-F^*_p;f_{p,sl}-f^*_p)\|, \quad
\omega_{p,sl}=\min\{\epsilon^x_{p,sl},\epsilon^f_{p,sl}\}.$$
Here, $(x_{p,sl}, y_{p,sl})$ denotes the solution obtained by solver $sl$ on problem $p$ with upper and lower level objective values $F_{p,sl}, f_{p,sl}$, $(x_p^*,y_p^*)$ is the known global optimal solution with objective values $F^*_p, f^*_p$. 
The performance ratio with respect to $\omega_{p,sl}$ is defined as:
\begin{equation*}
	r_{p,sl} =\begin{case}
		\log_{10}\left(\frac{\omega_{p,sl}}{\min_{sl'\in\mathcal{SL}}\omega_{p,sl'}}\right), &\text{if } \min_{sl'\in\mathcal{SL}}\omega_{p,sl'}\ne0 \text{ and } 0<\omega_{p,sl}\le10^{-3},\\
		0, &\text{if } \min_{sl'\in\mathcal{SL}}\omega_{p,sl'}=0 \text{ and }\omega_{p,sl}=0 ,\\
		50, &\text{if } \min_{sl'\in\mathcal{SL}}\omega_{p,sl'}=0\text{ and } 0<\omega_{p,sl}\le10^{-3},\\
		100,&\text{otherwise}.
	\end{case}
\end{equation*}
Similarly, we define the counterpart for CPU time as $r_{p,sl} =\log_{10}\left(\frac{t_{p,sl}}{\min_{sl'\in\mathcal{SL}}t_{p,sl'}}\right)$, where $t_{p,sl}$ represents the CPU time obtained by solver $sl$ for problem $p$.
The cumulative performance profile is defined by:
\begin{equation*}\label{p_sl}
	P_{sl}(\gamma) = \frac{|\{p \in \mathcal{P}: r_{p,sl} \leq \gamma\}|}{n_p}, \quad \gamma \geq 0.
\end{equation*}

   \subsection{Results and Analysis}        
Through our experiments, improved reference solutions were obtained for eight BOLIB problems:
AiyoshiShimizu1984Ex2 ($x^*=[0;0],y^*=[-10;-10]$), 
Dempe1992a ($x^*=[0;0],y^*=[0;0]$), 
Dempe1992b ($x^*=0,y^*=0$), 
DempeDutta2012Ex31 ($x^*=[1;1],y^*=[0;2]$), 
FalkLiu1995 ($x^*=[0.75;0.75],y^*=[0.75;0.75]$), 
MitsosBarton2006Ex318 ($x^*=1,y^*=0$), 
MorganPatrone2006c ($x^*=1.75,y^*=1$), 
ShimizuEtal1997a ($x^*=0.9459,y^*=-0.1622$).
All of these solutions strictly satisfy the original constraints and improve upon the reference points in the BOLIB library.

%

\begin{table}[!h]
	\centering
	\caption{Performance comparison of four algorithms }
	\begin{tabular}{|l|c|c|c|c|}
		\hline
		\textbf{Algorithm}  & ~~\textbf{Suc.}~~ & ~~\textbf{Obj.Suc.}~~ & ~~\textbf{Sol.Suc.}~~ &  ~~\textbf{Time (s)}~~  \\ \hline
		SVF-SBAL& \textbf{82} &\textbf{66} &\textbf{74}& 8.01e-02  \\ \hline
		KP-SBAL &  61&47& 55& 6.09e-02  \\ \hline
		KP-RLX & 64& 51& 61&{5.37e-02}  \\ \hline
		eRMDP &  70& 55& 65&\textbf{ 4.10e-02}  \\ \hline
	\end{tabular}
	\label{tab:algorithm_comparison}
\end{table}

		Table \ref{tab:algorithm_comparison} reports the performance of four competing methods over 109 benchmark problems. Evaluation is based on four criteria: 
		\textbf{Success (Sus.)}, requiring the composite error satisfying $\omega_{p,sl} \leq 10^{-3}$. Here, the composite error is defined as the minimum of the solution and objective errors, indicating that success is achieved when the solver attains satisfactory accuracy in either aspect; \textbf{Objective Success (Obj.Suc.)}, requiring the objective gap to be below tolerance: $\epsilon^f_{p,sl} \leq 10^{-3}$; \textbf{Solution Success (Sol.Suc.)}, requiring the solution error to be small: $\epsilon^x_{p,sl} \leq 10^{-3}$; and \textbf{Average Time (Time)} stands for the average CPU time exclusively over problems satisfying Success criteria.

\begin{figure}[!h]\centering
			\subfigure[Success Rate]{
				\includegraphics[width=7.3cm]{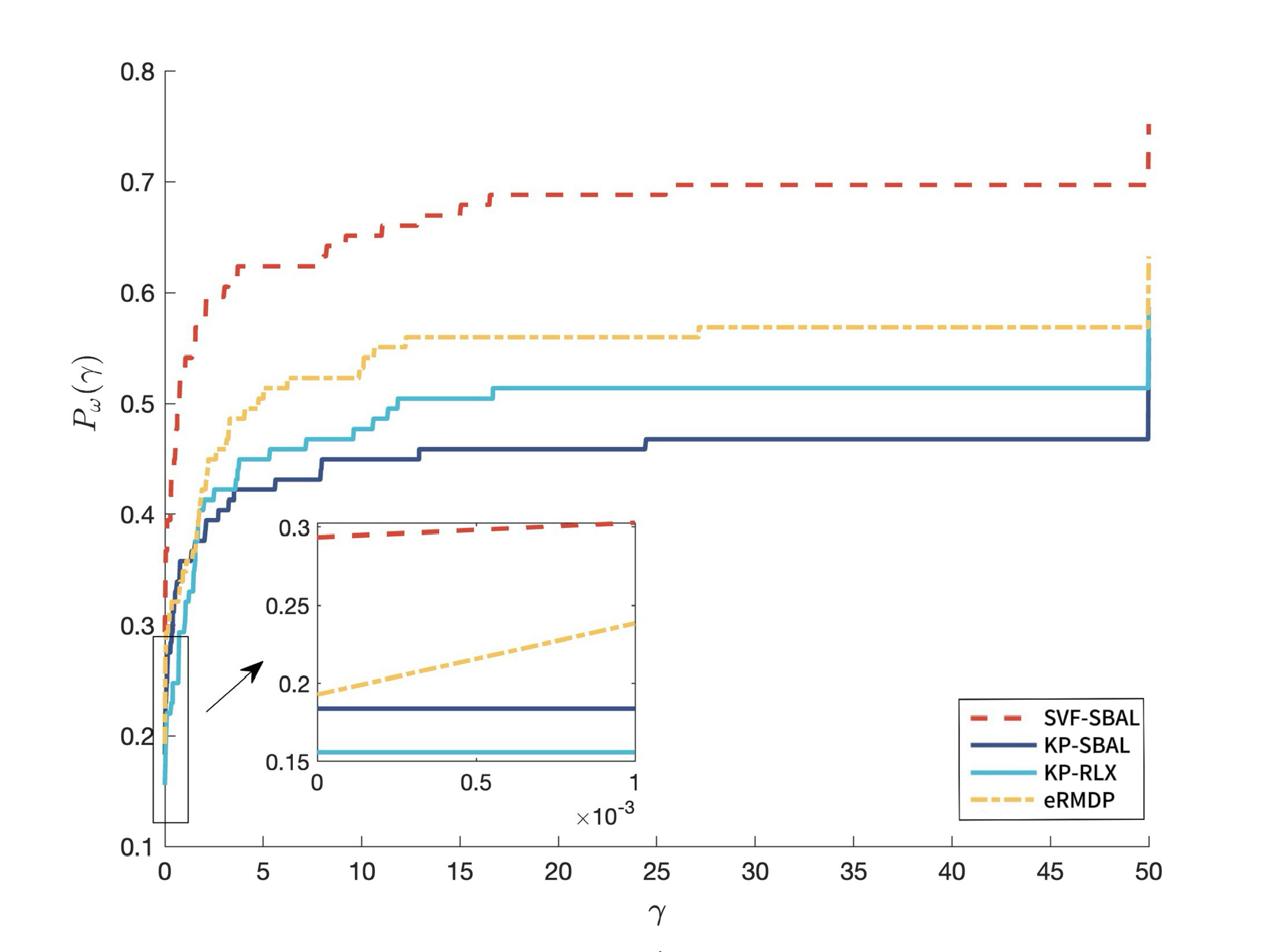}}\quad
			\subfigure[Objective function Rate]{
				\includegraphics[width=7.3cm]{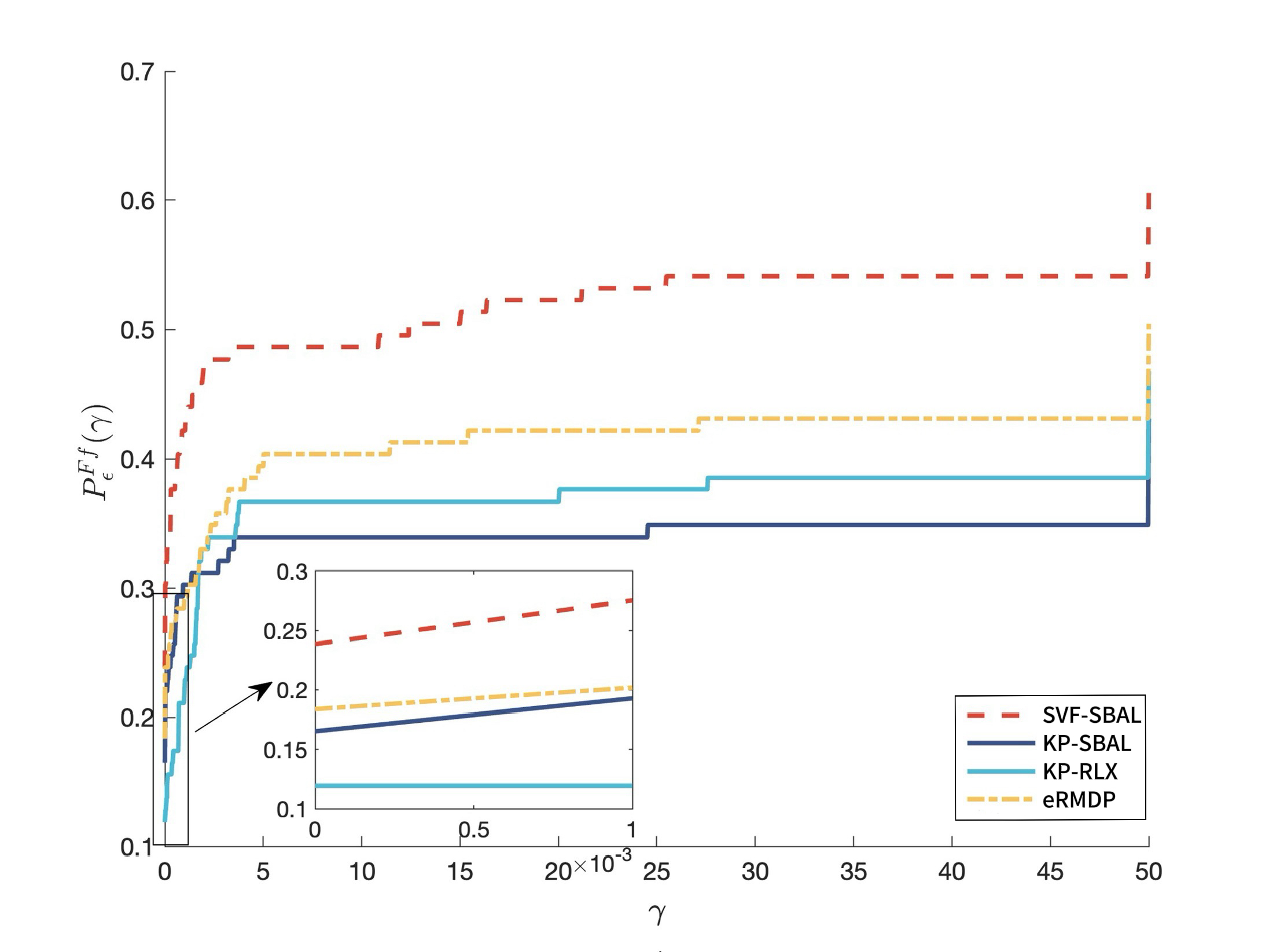}\label{fig:algorithm_performance3}
			}\quad
			\subfigure[Solution Rate]{
				\includegraphics[width=7.3cm]{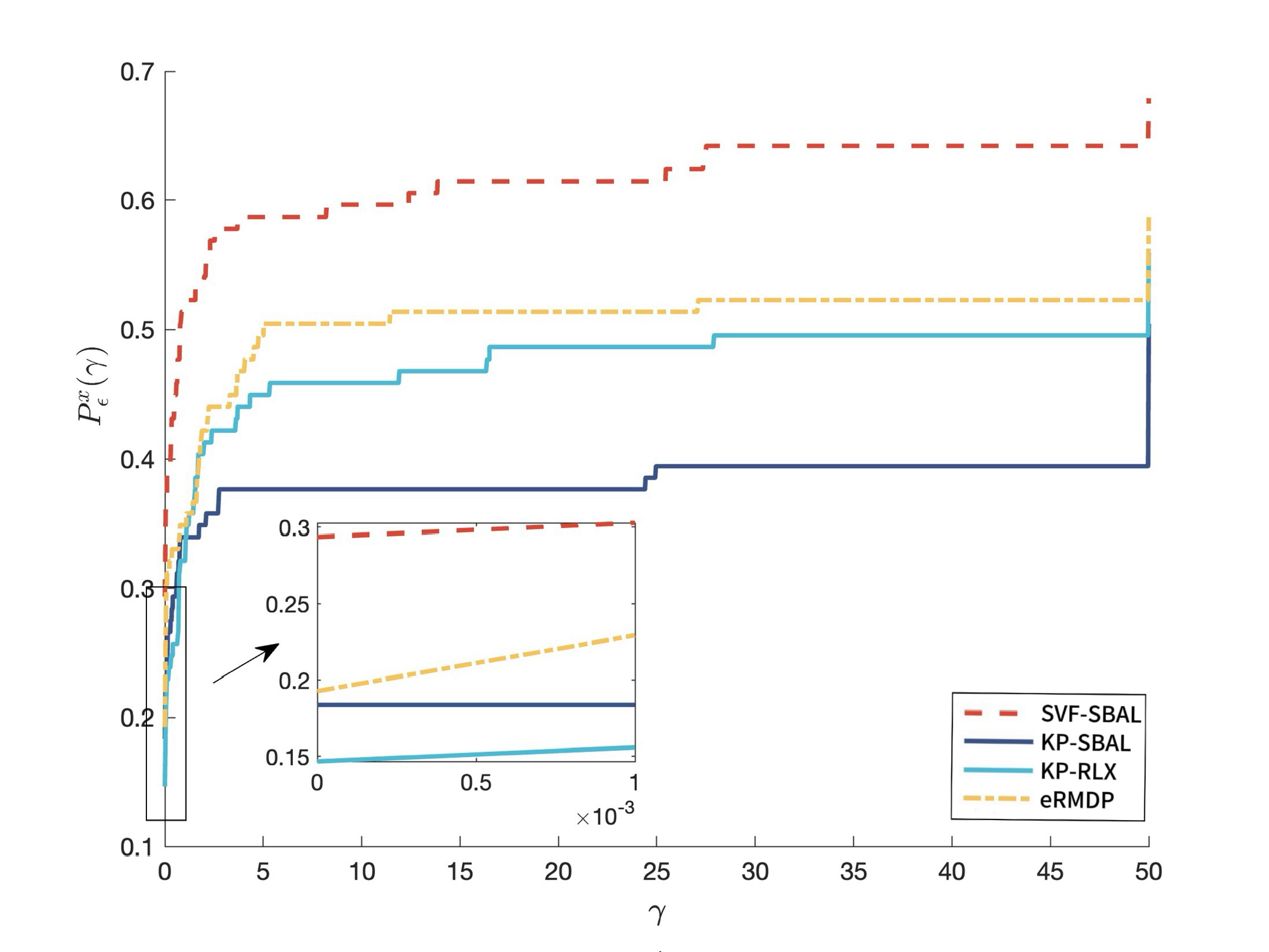}\label{fig:algorithm_performance2}
			}\quad
				\subfigure[Time Rate]{
				\includegraphics[width=7.3cm]{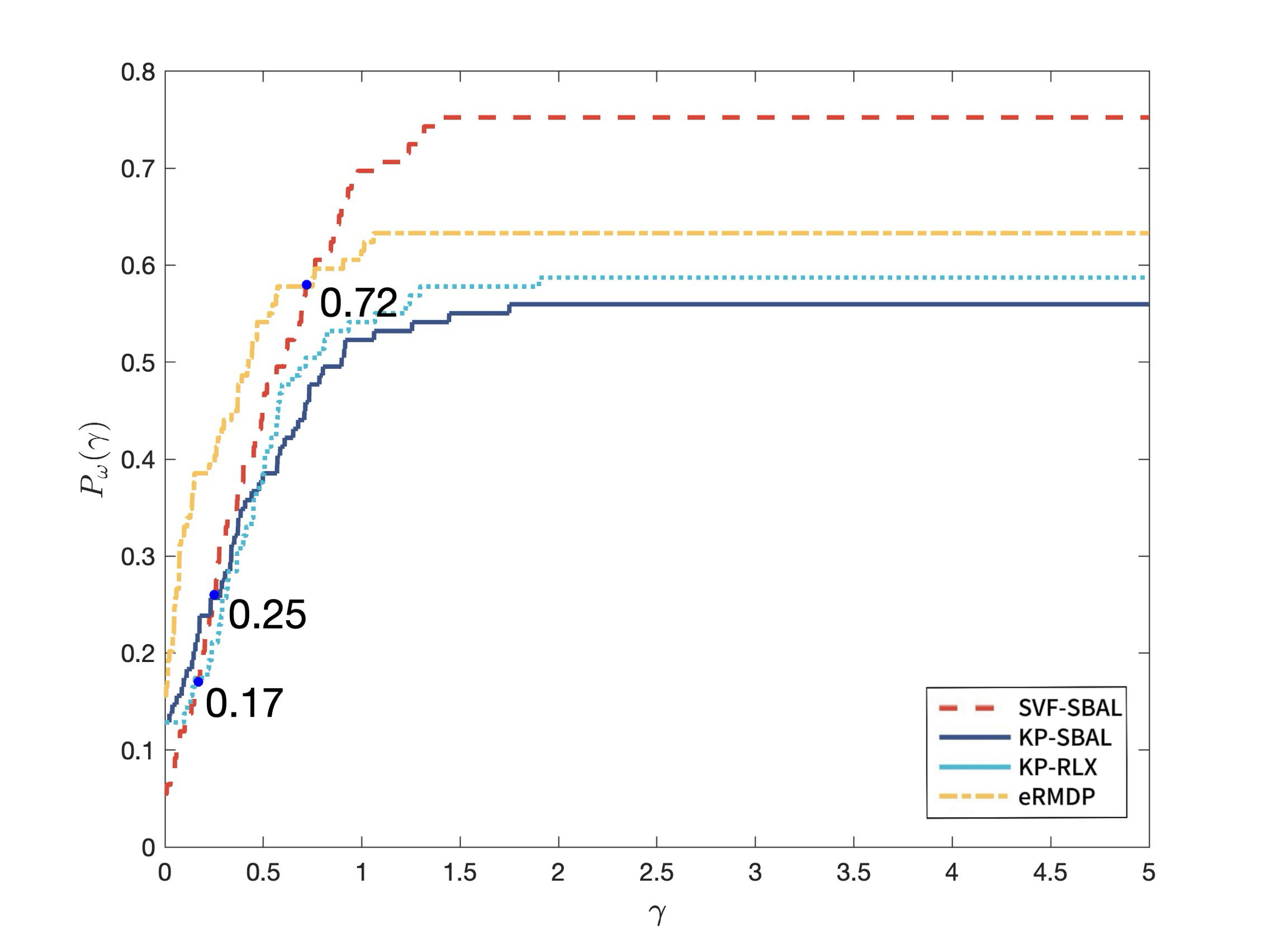}\label{fig:algorithm_performance4}
			}
			\caption{Performance profiles}
			\label{fig:algorithm_performance}
		\end{figure}

As shown in the table \ref{tab:algorithm_comparison}, SVF-SBAL attains the widest coverage and the best accuracy across all metrics: it achieves the largest counts for Success and leads on both objective- and solution-based criteria. 
The performance profiles in figure \ref{fig:algorithm_performance} reinforce this conclusion: the success curve for SVF-SBAL lies above all competitors and is the best performer on roughly 29\% of the instances; the objective-success curve and the solution-success curve show a similar lead and SVF-SBAL is best on about 24\% and 29\% of the problems, respectively.
The time-based performance profile in Fig. \ref{fig:algorithm_performance}(d) shows that SVF-SBAL solves more problems than KP-RLX, KP-SBAL, and eRMDP when allowed 1.48, 1.78, and 5.25 times the runtime of the fastest solver, respectively.

The model SVF is noticeably more robust in producing high-quality approximations on challenging bilevel instances, notably those with nonconvex lower levels (e.g. Mirrlees1999, Mitsos–Barton instances) or with violated lower level KKT conditions (e.g. Dempe1992b), where KKT model fails.
The \textbf{Mirrlees1999} instance provides a telling example: this classical nonconvex benchmark is known to be inaccessible to KKT formulation \cite{yz2}. On this instance, (SVF-SBAL) computes a high-accuracy approximate solution (solution error $\approx 1.54e-02$, objective error  $\approx 9.18e-04$), which illustrates that the method can recover near-optimal solutions in cases where standard KKT formulations fails.


In summary, these results demonstrate that the proposed SVF formulation, solved via the SBAL algorithm, achieves the highest overall success rates, superior accuracy, and competitive runtimes. More importantly, it solves benchmark problems where classical KKT formulation fundamentally fails.

\section{Conclusion}
We introduced the SVF formulation that replaces the implicit value function with an explicit surrogate defined at a lower level stationary point. 
This model preserves the comparison mechanism of VP while avoiding its intractability and the failure of KKT formulation.
Using the SBAL approximation, we established convergence of the smoothed problems to SVF solutions and Clarke stationary points. 
Numerical results on 109 benchmarks confirm that SVF consistently outperforms existing KP and MDP reformulations, including success on the classical Mirrlees' problem where KKT fomulation fail.


\end{document}